\documentclass[11pt]{amsart}

\usepackage{amssymb,amsthm,amsmath}
\usepackage[numbers,sort&compress]{natbib}
\usepackage{color}

\usepackage{graphicx}
\usepackage{tikz}
\usepackage[colorlinks,linkcolor=blue,anchorcolor=blue,citecolor=blue]{hyperref}

\newcommand{\SL}{\mathrm{SL}(2,\mathbb{R})}

\let\newpf\proof \let\proof\relax 
\newenvironment{pf}{\newpf[\proofname]}{\qed\endtrivlist}

\newcommand{\ba}{\overline{A}}

\def\be{\begin{equation}}
\def\ee{\end{equation}}

\def\ba{{\begin{align}}}
\def\ea{{\end{align}}}

\def\bm{\begin{matrix}}
\def\em{\end{matrix}}

\def\SL{{\mathrm{SL}}}

\def\0{{\mathbf 0}}

\newtheorem{Theorem}{Theorem}[section]
\newtheorem*{Theorem*}{Theorem 1.1}
\newtheorem{Lemma}[Theorem]{Lemma}
\newtheorem{Proposition}[Theorem]{Proposition}
\newtheorem{Corollary}[Theorem]{Corollary}

\newtheorem{Claim}{Claim}

\theoremstyle{definition}
\newtheorem{Remark}[Theorem]{Remark}

\numberwithin{equation}{section}

\theoremstyle{definition}

\renewcommand{\mod}{\operatorname{mod}}

\newcommand{\supp}{\operatorname{supp}}
\newcommand{\id}{\operatorname{id}}

\newcommand{\C}{{\mathbb C}}
\newcommand{\D}{{\mathbb D}}

\newcommand{\Q}{{\mathbb Q}}
\newcommand{\R}{{\mathbb R}}
\newcommand{\T}{{\mathbb T}}

\newcommand{\Z}{{\mathbb Z}}

\def\B0{{\bold{0}}}

\catcode`\@=12

\def\Empty{}
\newcommand\oplabel[1]{
  \def\OpArg{#1} \ifx \OpArg\Empty {} \else
    \label{#1}
  \fi}

\newcommand{\comm}[1]{}
\newcommand{\comment}[1]{}

\begin{document}

\title[A.C.\ Spectrum for Quasi-Periodic CMV Matrices]{Absolutely Continuous Spectrum for CMV Matrices With Small Quasi-Periodic Verblunsky Coefficients}

\author{Long Li} \address{
Department of Mathematics, Nanjing University, Nanjing 210093, China
}

\email{huanzhensu@icloud.com}

\author{David Damanik}
\address{ Department of Mathematics, Rice University, Houston, Texas, 77005
}
\email{damanik@rice.edu}
\thanks{D.D.\ was supported in part by NSF grant DMS--1700131, an Alexander von Humboldt Foundation research award, and a Simons Fellowship}

\author{Qi Zhou}
\address{
Chern Institute of Mathematics and LPMC, Nankai University, Tianjin 300071, China
}

 \email{qizhou@nankai.edu.cn}
\thanks{
Q.Z.\ was  supported by
 National Key R\&D Program of China (2020YFA0713300), NSFC grant (12071232), The Science Fund for Distinguished Young Scholars of Tianjin (No. 19JCJQJC61300) and Nankai Zhide Foundation.}

\setcounter{tocdepth}{1}

\begin{abstract}
We consider standard and extended CMV matrices with small quasi-periodic Verblunsky coefficients and show that on their essential spectrum, all spectral measures are purely absolutely continuous. This answers a question of Barry Simon from 2005.
\end{abstract}

\maketitle

\section{Introduction}

This paper is concerned with the spectral analysis of (standard and extended) CMV matrices with quasi-periodic Verblunsky coefficients. CMV matrices are canonical matrix representations of unitary operators with a cyclic vector, and they arise naturally in the context of orthogonal polynomials on the unit circle. We refer the reader to \cite{Simon1, Simon2} for background.

Let us recall how CMV matrices arise in connection with orthogonal polynomials on the unit circle. Suppose $\mu$ is a non-trivial probability measure on the unit circle $\partial \mathbb{D} = \{ z \in \mathbb{C} : |z| = 1 \}$, that is, $\mu(\partial \D) = 1$ and $\mu$ is not supported on a finite set. By the non-triviality assumption, the functions $1$, $z$, $z^2,\cdots$ are linearly independent in the Hilbert space $\mathcal{H} = L^2(\partial\mathbb{D}, d\mu)$, and hence one can form, by the Gram-Schmidt procedure, the \textit{monic orthogonal polynomials} $\Phi_n(z)$, whose \textit{Szeg\H{o} dual} is defined by $\Phi_n^{*} = z^n\overline{\Phi_n({1}/{\overline{z}})}$. There are constants $\{\alpha_n\}_{n=0}^\infty$ in $\mathbb{D}=\{z\in\mathbb{C}:|z|<1\}$, called the \textit{Verblunsky coefficients}, so that
\begin{equation}\label{eq1}
\Phi_{n+1}(z) = z \Phi_n(z) - \overline{\alpha}_n \Phi_n^*(z),
\end{equation}
which is the so-called \textit{Szeg\H{o} recurrence}. Conversely, every sequence $\{\alpha_n\}_{n=0}^\infty$ in $\mathbb{D}$ arises in this way.

The orthogonal polynomials may or may not form a basis of $\mathcal{H}$. However, if we apply the Gram-Schmidt procedure to $1, z, z^{-1}, z^2, z^{-2}, \ldots$, we will obtain a basis -- called the \textit{CMV basis}. In this basis, multiplication by the independent variable $z$ in $\mathcal{H}$ has the matrix representation
\begin{equation*}
\mathcal{C}=\left(
\begin{matrix}
\overline{\alpha}_0&\overline{\alpha}_1\rho_{0}&\rho_1\rho_0&0&0&\cdots&\\
\rho_0&-\overline{\alpha}_1\alpha_{0}&-\rho_1\alpha_0&0&0&\cdots&\\
0&\overline{\alpha}_2\rho_{1}&-\overline{\alpha}_2\alpha_{1}&\overline{\alpha}_3\rho_2&\rho_3\rho_2&\cdots&\\
0&\rho_2\rho_{1}&-\rho_2\alpha_{1}&-\overline{\alpha}_3\alpha_2&-\rho_3\alpha_2&\cdots&\\
0&0&0&\overline{\alpha}_4\rho_3&-\overline{\alpha}_4\alpha_3&\cdots&\\
\cdots&\cdots&\cdots&\cdots&\cdots&\cdots&
\end{matrix}
\right),
\end{equation*}
where
\begin{equation}\label{e.rhodef}
\rho_n = (1-|\alpha_n|^2)^{1/2}
\end{equation}
for $n \ge 0$. A matrix of this form is called a \textit{CMV matrix}.

It is sometimes helpful to also consider a two-sided extension of a matrix of this form. Namely, given a bi-infinite sequence $\{ \alpha_n \}_{n \in \Z}$ in $\D$ (and defining the $\rho_n$'s as before), we may consider the \textit{extended CMV matrix}
\begin{equation*}
\mathcal{E}=\left(
\begin{matrix}
\cdots&\cdots&\cdots&\cdots&\cdots&\cdots&\cdots\\
\cdots&-\overline{\alpha}_0\alpha_{-1}&\overline{\alpha}_1\rho_{0}&\rho_1\rho_0&0&0&\cdots&\\
\cdots&-\rho_0\alpha_{-1}&-\overline{\alpha}_1\alpha_{0}&-\rho_1\alpha_0&0&0&\cdots&\\
\cdots&0&\overline{\alpha}_2\rho_{1}&-\overline{\alpha}_2\alpha_{1}&\overline{\alpha}_3\rho_2&\rho_3\rho_2&\cdots&\\
\cdots&0&\rho_2\rho_{1}&-\rho_2\alpha_{1}&-\overline{\alpha}_3\alpha_2&-\rho_3\alpha_2&\cdots&\\
\cdots&0&0&0&\overline{\alpha}_4\rho_3&-\overline{\alpha}_4\alpha_3&\cdots&\\
\cdots&\cdots&\cdots&\cdots&\cdots&\cdots&\cdots
\end{matrix}\right).
\end{equation*}

Naturally, one is interested in both direct and inverse spectral results, depending on whether one starts with information about the Verblunsky coefficients or the measure. This paper is concerned with a direct spectral problem. The Verblunsky coefficients will be small and quasi-periodic, and our goal is to show that the associated spectral measures are purely absolutely continuous (on the essential spectrum). We aim to establish this property for both standard and extended CMV matrices. This addresses one of the open problems described by Simon in \cite{Simon2}. Indeed, in the Remarks and Historical Notes to \cite[Section~10.16]{Simon2}, he writes that from his discussion of ergodic Verblunsky coefficients ``conspicuously absent is the case of almost periodic Verblunsky coefficients'' and ``especially interesting is the quasiperiodic case.'' He goes on to suggest that one should prove absolute continuity for small quasi-periodic coefficients and pure point spectrum for some quasi-periodic examples. We note that quasi-periodic examples with pure point spectrum were exhibited by Wang-Damanik in \cite{DW}, making use of earlier results of Zhang \cite{Zhang1} on the positivity of the Lyapunov exponent in the setting in question. Here we address the first part of Simon's question, namely how to prove absolute continuity of the spectral measures for small quasi-periodic Verblunsky coefficients.

Let us now describe the setting and the main result in detail. We consider small analytic quasi-periodic Verblunsky coefficients of the form
\begin{equation}\label{setting}
\alpha_{n}(x) = \alpha(x+(n-1)\omega), \quad n\in \mathbb{Z},
\end{equation}
where
\begin{equation}\label{setting2}
\alpha(x) = \lambda e^{2\pi i h(x)},
\end{equation}
$h \in C^{\omega}(\mathbb{T}^{d},\mathbb{R})$, $\lambda \in(0,1)$, $\omega \in\mathbb{R}^{d}$ with $\langle m,\omega\rangle\notin \mathbb{Z}$ for any $m\in \mathbb{Z}^{d}\backslash\{0\}$.

The sequence in \eqref{setting}--\eqref{setting2} defines an extended CMV matrix. Since this matrix formally depends on $x \in \mathbb{T}^{d}$, we will denote it by $\mathcal{E}_x$. While it of course also depends on $h$ and $\omega$, we view them as fixed and suppress them from the notation. To define a standard CMV matrix, we only consider the values of $\alpha_{n}(x)$ in \eqref{setting} for $n \in\mathbb{Z}_{+}$ and denote the resulting one-sided CMV matrix by $\mathcal{C}_x$.

Recall that $\omega\in\mathbb{R}^{d}$ is called {\it Diophantine} if there exist some $\kappa, \tau > 0$ such that
$$
\inf \limits_{j \in\mathbb{Z}} \vert \langle n,\omega\rangle - j \vert \geq \frac{\kappa}{\vert n\vert^{\tau}}
$$
for all $n \in \mathbb{Z}^{d} \backslash \{0\}$. Let $\mathrm{DC}(\kappa,\tau)$ be the set of all Diophantine numbers with prescribed $\kappa,\tau$.

We can now state our main result:

\begin{Theorem}\label{THM4.3}
Suppose that $r, \kappa, \tau > 0$, $h \in C^{\omega}_r(\mathbb{T}^{d},\mathbb{R})$, and $\omega \in \mathrm{DC}(\kappa,\tau)$. Then there exists $\lambda_{0}=\lambda_{0}(r, \kappa, \tau) > 0$ such that for $\lambda \in (0, \lambda_{0})$, the following holds for the Verblunsky coefficients given by \eqref{setting}--\eqref{setting2}:
\begin{itemize}

\item For every $x \in \mathbb{T}^{d}$, $\mathcal{E}_x$ has purely absolutely continuous spectrum.

\item For every $x \in \mathbb{T}^{d}$, the restriction of the canonical spectral measure of $\mathcal{C}_x$ to $\Sigma = \sigma_\mathrm{ess}(\mathcal{C}_x)$ is purely absolutely continuous and $\mathcal{C}_{x}$ has at most one eigenvalue in each connected component of $\partial\mathbb{D}\backslash \Sigma$.

\end{itemize}
\end{Theorem}

\begin{Remark}
(a) The set of all Diophantine $\omega \in \mathbb{R}^{d}$ has full Lebesgue measure. It is an interesting question whether the conclusion of the theorem can fail for some $\omega$'s, even if $\lambda$ is small. In the Schr\"odinger setting, it is known that the cases $d = 1$ and $d > 1$ behave differently: there is a non-perturbative version of this result (which means that $\lambda_{0} > 0$ can be chosen uniformly in $\kappa$ and $\tau$ and the result then holds for all Diophantine $\omega$ with this choice of $\lambda_0$) in the case $d = 1$ \cite{AJ1}, and it is known that there cannot be a non-perturbative version of this result when $d > 1$ \cite{B02}. It would be interesting to establish results of this kind in the OPUC setting as well. We explain in Appendix~\ref{app.NP} how one can approach the first problem.

(b) While we state the result in the analytic topology, our approach works in the $C^{\infty}$ topology since our proof is KAM based.\footnote{Fayad-Krikorian \cite{FK} provide the $C^\infty$ version of the KAM result underlying our approach.} Thus, we exactly solve the question Simon proposed in the Remarks and Historical Notes to \cite[Section~10.16]{Simon2}.

(c) Going beyond the explicit formulation of Simon's question, the curious reader may wonder whether in addition to the pure point result from \cite{DW} and the absolute continuity result in Theorem~\ref{THM4.3} one can also exhibit analytic (or smooth) quasi-periodic CMV matrices with purely singular continuous spectrum. The answer is yes and all necessary tools to produce such examples (of \textit{extended} analytic quasi-periodic CMV matrices with purely singular continuous spectrum) already exist. We include a discussion to this effect in Appendix~\ref{app.sc} for the sake of said curious reader.
\end{Remark}

We remark that from the OPUC perspective, standard CMV matrices are the natural object and hence we are primarily interested in the spectral properties of the matrices $\mathcal{C}_x$. However, from the perspective of the actual spectral analysis of these matrices, it is crucial to also consider the extended CMV matrices $\mathcal{E}_x$ since the underlying torus translation
$$
T_\omega : \mathbb{T}^{d} \to \mathbb{T}^{d}, \quad x \mapsto x + \omega
$$
is an \textit{invertible ergodic transformation}, which enables us to rely on the powerful general theory that has been developed for extended CMV matrices with Verblunsky coefficients generated by continuous sampling along the orbit of such a transformation; compare \cite{DFLY, FDG} and references therein. We will be more specific below, but merely mention here that the spectrum can be characterized in terms of the uniform hyperbolicity of the associated Szeg\H{o} cocycles \cite{DFLY} and the absolutely continuous spectrum can be investigated via Kotani theory \cite{FDG, FO, G93, GT94, Kotani, Simon2} (and both of these fundamental tools require the two-sided setting in their full-fledged version). Consequently, we will prove the desired absolute continuity of spectral measures both for the standard CMV matrices $\mathcal{C}_x$ and for the extended CMV matrices $\mathcal{E}_x$. To facilitate stating our results, let us mention now that given $h$, $\lambda$, and $\omega$, there is a compact set $\Sigma \subseteq \partial \mathbb{D}$ such that for every $x \in \mathbb{T}^{d}$, we have
\begin{equation}\label{e.spectrum}
\Sigma = \sigma_\mathrm{ess}(\mathcal{C}_x) = \sigma_\mathrm{ess}(\mathcal{E}_x) = \sigma(\mathcal{E}_x).
\end{equation}
The discrete spectrum of $\mathcal{C}_x$ may (and usually will) be $x$-dependent, but it is known that each connected component of $\partial \mathbb{D} \backslash \Sigma$ contains at most one point that belongs to $\sigma_\mathrm{disc}(\mathcal{C}_x)$. We will say more about this later in the paper.

It is well known, and crucial to our present discussion, that the theory of orthogonal polynomials on the unit circle shares many similarities with the theory of one-dimensional discrete Schr\"odinger operators. This is not a coincidence as the latter embeds naturally into the theory of orthogonal polynomials on the real line, and the unit circle and the real line stand out as those subsets of the complex plane on which polynomials that are orthogonal with respect to a non-trivial probability measure that is supported on one of these two sets are known to obey a useful recursion (the Verblunsky coefficients are just the recursion coefficients in the case of the unit circle; compare \eqref{eq1}). As a result, there has been an extensive effort to work out OPUC analogs of results that had first been established for discrete Schr\"odinger operators. From this perspective, the result of Zhang mentioned above \cite{Zhang1} is an OPUC analog of results of Herman \cite{H} and Sorets-Spencer \cite{SS91} in the Schr\"odinger setting, and the Wang-Damanik result \cite{DW} is an analog of a result of Bourgain-Goldstein \cite{BG}. The results of the present paper are the OPUC analog of a Schr\"odinger operator result that has evolved over the years with many important contributions. We will give more detailed pointers to the literature below, but mention here some of the milestones: Dinaburg-Sinai \cite{DS}, Eliasson \cite{Eli92}, Avila-Jitomirskaya \cite{AJ1},  Avila \cite{Avi}.

We wish to emphasize that working out the OPUC analog of a Schr\"odinger result is not straightforward. There is no general mechanism that serves as a black box transforming an input into an output without further work. Rather, the plethora of known companion results exists primarily because many of the tools of spectral analysis on the Schr\"odinger side have been shown to have analogs on the OPUC side, and using the latter allowed the authors in question to obtain the desired OPUC results by following a similar proof strategy.

In the case of the problem at hand, a crucial tool in the spectral analysis of the quasi-periodic Schr\"odinger case, namely Aubry duality, is not known to exist in a similarly useful form in the OPUC setting. Aubry duality allows one, in the Schr\"odinger case, to relate the absolute continuity problem via duality to a localization problem for the dual model and to exploit the powerful methods that have been developed to establish such localization statements; we refer the reader to  \cite{Avi,AJ1,BJ1}  for results of this kind.
The existence of absolutely continuous spectrum is well known to follow from the boundedness of the solutions of the difference equation associated with the operator in question, and this boundedness can often be established in the quasi-periodic setting by describing these solutions with the help of suitable cocycles, and then conjugating these cocycles to $\mathrm{SO}(2,\R)$-valued cocycles. The key issue is the pure absolute continuity as one typically obtains boundedness only on a set of spectral parameters that has full Lebesgue measure, and hence the absence of singular spectrum still needs to be clarified. To this end we follow here an approach initially developed by Avila \cite{Avi}. In this approach, one needs to prove that the set of spectral parameters with unbounded solutions has zero weight with respect to the spectral measure, which in turn relies on a measure estimate for the set of spectral parameters for which the corresponding cocycles have a given growth rate. In particular, one seeks to establish suitable almost reducibility results (where one can conjugate into an arbitrarily small neighborhood of a constant). The duality approach shows that (almost) reducibility for the initial model can be derived from (almost) localization for the dual model.   As we do not have the duality approach at our disposal for the problem we study, we have to look for direct methods of proving absolute continuity via (almost) reducibility.  Instead we establish a quantitative almost reducibility result directly with the help of recent advances by Cai-Chavaudret-You-Zhou \cite{CCYZ} and Leguil-You-Zhao-Zhou \cite{LYZZ}\footnote{This approach was first used to deal with the absolutely continuous spectrum of Schr\"odinger operators with quasi-periodic-like potentials \cite{XYZ}. }. We also mention that there is another proof of purely absolutely continuous spectrum for continuum Schr\"odinger operators with small analytic quasi-periodic potentials by Eliasson \cite{Eli92} that does not use duality. Compared with \cite{Eli92} our approach is more concise: we do not need to complexify the energy and estimate the imaginary part of the Green function $G(E\pm i g; n,m)$ (which corresponds to $G((1\pm \delta)e^{i\zeta}; n,m)$ in our case).

Finally, we mention that in the Schr\"odinger context, the approach to the study of the spectral properties of a quasi-periodic operator that is based on quantitative almost reducibility is very fruitful \cite{Avi,AJ1,ayz,CCYZ,Eli92,LYZZ} (see also the nice survey of You \cite{You} for more results). Our result is the first realization of this in the OPUC setting.

The structure of the paper is as follows. After recalling some relevant parts of the general OPUC theory in Section~\ref{sec.2}, we prove a quantitative almost reducibility result for analytic quasi-periodic $\mathrm{SU}(1,1)$ cocycles (derived from the work of Cai-Chavaudret-You-Zhou \cite{CCYZ} and Leguil-You-Zhao-Zhou \cite{LYZZ}, which we briefly recall in Appendix~\ref{sec.app}) in Section~\ref{sec.3} and a lower bound for the density of states measure in the quasi-periodic case we study in Section~\ref{sec.4}. The proof of Theorem~\ref{THM4.3} is given in Section~\ref{sec.5}. Finally, we discuss how to obtain a non-perturbative result in the case $d = 1$ in Appendix~\ref{app.NP} and how to obtain analytic quasi-periodic extended CMV matrices with singular continuous spectrum in Appendix~\ref{app.sc}.

\section{Preliminaries}\label{sec.2}

In this section we collect some material we will need in the subsequent sections. The results we describe here are well known, but they are included for the convenience of the reader.

\subsection{The Standard Factorization of CMV Matrices}

Recall that, given a sequence $\{\alpha_{n}\}_{n\in\mathbb{Z}_{+}}\subset\mathbb{D}=\{z\in\mathbb{C}:|z|<1\}$, the standard (or half-line) CMV matrix takes the form
$$
\mathcal{C}=\left(
\begin{matrix}
\overline{\alpha}_0&\overline{\alpha}_1\rho_{0}&\rho_1\rho_0&0&0&\cdots&\\
\rho_0&-\overline{\alpha}_1\alpha_{0}&-\rho_1\alpha_0&0&0&\cdots&\\
0&\overline{\alpha}_2\rho_{1}&-\overline{\alpha}_2\alpha_{1}&\overline{\alpha}_3\rho_2&\rho_3\rho_2&\cdots&\\
0&\rho_2\rho_{1}&-\rho_2\alpha_{1}&-\overline{\alpha}_3\alpha_2&-\rho_3\alpha_2&\cdots&\\
0&0&0&\overline{\alpha}_4\rho_3&-\overline{\alpha}_4\alpha_3&\cdots&\\
\cdots&\cdots&\cdots&\cdots&\cdots&\cdots&
\end{matrix}
\right),
$$
where $\rho_{j}=\sqrt{1-\vert \alpha_{j}\vert^{2}}$. It is possible to factorize this matrix as follows,
$$
\mathcal{C}=\mathcal{L}\mathcal{M},
$$
where
$$
\mathcal{L}=\begin{pmatrix}\Theta_{0}\\
&\Theta_{2}\\
&&\Theta_{4}\\
&&&\ddots\end{pmatrix},\quad \mathcal{M}=\begin{pmatrix}
1 &&\\
&\Theta_{1}&\\
&&\Theta_{3}\\
&&& \ddots
\end{pmatrix}
$$
and
$$
\Theta_{j}=\begin{pmatrix} \bar{\alpha}_{j} &\rho_{j}\\ \rho_{j}  &-\alpha_{j} \end{pmatrix}.
$$

Similarly, the extended CMV matrix
$$
\mathcal{E}=\left(
\begin{matrix}
\cdots&\cdots&\cdots&\cdots&\cdots&\cdots&\cdots\\
\cdots&-\overline{\alpha}_0\alpha_{-1}&\overline{\alpha}_1\rho_{0}&\rho_1\rho_0&0&0&\cdots&\\
\cdots&-\rho_0\alpha_{-1}&-\overline{\alpha}_1\alpha_{0}&-\rho_1\alpha_0&0&0&\cdots&\\
\cdots&0&\overline{\alpha}_2\rho_{1}&-\overline{\alpha}_2\alpha_{1}&\overline{\alpha}_3\rho_2&\rho_3\rho_2&\cdots&\\
\cdots&0&\rho_2\rho_{1}&-\rho_2\alpha_{1}&-\overline{\alpha}_3\alpha_2&-\rho_3\alpha_2&\cdots&\\
\cdots&0&0&0&\overline{\alpha}_4\rho_3&-\overline{\alpha}_4\alpha_3&\cdots&\\
\cdots&\cdots&\cdots&\cdots&\cdots&\cdots&\cdots
\end{matrix}\right)
$$
can be written as
$$
\mathcal{E} = \mathcal{L}'\mathcal{M}',
$$
where
$$
\mathcal{L}'=\begin{pmatrix}\ddots\\  &\Theta_{-2}\\ & &\Theta_{0}\\
 & & &\Theta_{2}\\
& & &&\ddots\end{pmatrix},\quad \mathcal{M}'=\begin{pmatrix}
\ddots \\
&\Theta_{-1} \\
 & &\Theta_{1}\\
&&& \ddots
\end{pmatrix}.
$$

\subsection{Maximal Spectral Measures}\label{MSM}

By construction, the orthogonality measure $\mu$ with respect to which the CMV matrix $\mathcal{C}$ is constructed is a maximal spectral measure for this matrix. In other words, for every $\psi \in \ell^2(\Z_+)$, the spectral measure corresponding to the pair $(\mathcal{C},\psi)$ is absolutely continuous with respect to $\mu$. In particular, $\mathcal{C}$ has spectral multiplicity one. An extended CMV matrix $\mathcal{E}$, on the other hand, in general does not have spectral multiplicity one. However, the spectral multiplicity of $\mathcal{E}$ is at most two and the spectral subspace corresponding to $\mathcal{E}$ and the pair of vectors $\delta_0, \delta_1$ is all of $\ell^2(\Z)$ \cite[Theorem~10.16.5]{Simon2}. Thus, the sum of the spectral measures corresponding to $(\mathcal{E},\delta_0)$ and $(\mathcal{E},\delta_1)$ is maximal and will be denoted by $\Lambda$.

For our later discussion let us keep in mind that in order to prove that all spectral measures of a standard (resp., extended) CMV matrix $\mathcal{C}$ (resp., $\mathcal{E}$) are purely absolutely continuous on some set $S \subseteq \partial \mathbb{D}$, it suffices to show that $\mu$ (resp., $\Lambda$) is purely absolutely continuous on $S$.

Note also that for the quasi-periodic CMV matrices $\mathcal{C}_x$ and $\mathcal{E}_x$ we study, the maximal spectral measures of course depend on $x \in \mathbb{T}^d$ and will be denoted by $\mu_x$ and $\Lambda_x$, respectively.

\subsection{Transfer Matrices}

By normalizing the monic orthogonal polynomials $\Phi_n(z)$ in $L^2(\partial\mathbb{D}, d\mu)$, we obtain the orthonormal polynomials $\varphi_n(z)$. Their Szeg\H{o} duals $\varphi_n^*(z)$ are defined as before. The recursion \eqref{eq1} can then be rephrased as follows,
\begin{equation}\label{e.szegorecnorm}
\begin{pmatrix} \varphi_{n}(z) \\ \varphi^{*}_{n} (z)\end{pmatrix} = \frac{1}{\rho_{n}} \begin{pmatrix} z & -\bar{\alpha}_{n} \\ -\alpha_{n}z & 1 \end{pmatrix} \begin{pmatrix} \varphi_{n-1}(z) \\ \varphi^{*}_{n-1}(z) \end{pmatrix} , \varphi_{0}(z)=\varphi^{*}_{0}(z)=1,
\end{equation}
where the $\alpha_n$ are the Verblunsky coefficients and $\rho_n$ is given by \eqref{e.rhodef}.

The matrix
$$
\tilde{S}(\alpha,z) = \frac{1}{\sqrt{1-\vert\alpha\vert^{2}}} \begin{pmatrix} z&-\bar{\alpha} \\ -\alpha z & 1 \end{pmatrix},
$$
which provides one step on the iteration of \eqref{e.szegorecnorm}, is the called the \textit{Szeg\H{o} cocycle map}.

In this paper, we find it convenient to work with the \textit{renormalized Szeg\H{o} cocycle map}
\begin{equation}\label{e.renszegococyclemap}
S(\alpha,z)=z^{-\frac{1}{2}}\tilde{S}(\alpha,z) \in \mathrm{SU}(1,1).
\end{equation}

We need the following relations between the solutions of the generalized eigenvalue equation of $\mathcal{C}_{x}$ and the Szeg\H{o} polynomials.
 Suppose $u, v$ are defined as $$
\mathcal{C}_{x}u=zu,\quad v=\mathcal{M}u,
$$
that is, $u$ is the solution of the generalized eigenvalue equation of $\mathcal{C}_{x}$ and $v$ is obtained by applying a unitary transformation on $u$.  The vector $\begin{pmatrix}u_{n}\\ v_{n}\end{pmatrix}$ then obeys the Gesztesy-Zinchenko iterations. It relates to the Szeg\H{o} cocycle via the following:
\begin{equation}\label{e.sol1}
u_{2n}=z^{-n}\varphi_{2n}(z),u_{2n+1}=z^{-(n+1)}\varphi^{*}_{2n+1}(z)
\end{equation}
\begin{equation}\label{e.sol2}
v_{2n}=z^{-n}\varphi^{*}_{2n}(z),v_{2n+1}=z^{-n}\varphi_{2n+1}(z)
\end{equation}
with $\begin{pmatrix}u_{0}\\v_{0}\end{pmatrix}=\begin{pmatrix}1\\1\end{pmatrix}$ as boundary values, for all $z \in \mathbb{C}\backslash \{0\}$ and $n\geq 0$ (see \cite[(3.4)--(3.7)]{GDO}). Specifically, for $z \in \partial\mathbb{D},$ we have $\vert\varphi_{n}\vert=\vert \varphi^{*}_{n}\vert,$ which implies
\begin{equation}\label{equivalentmodel}
\vert u_{n}\vert=\vert \varphi_{n}\vert, n\in\mathbb{N}.
\end{equation}
For extended CMV matrices, there are very similar results. Let $s,t$ be the solutions of the following generalized eigenvalue equations
$$
\mathcal{E}_{x}s = zs, \quad t=\mathcal{M}'s.
$$
Then $\begin{pmatrix}s_{n}\\t_{n}\end{pmatrix}$ obeys the Gesztesy-Zinchenko iteration: for $n \in \mathbb{Z}$,
\begin{equation}\label{e.GZ1}\begin{pmatrix}s_{n}\\ t_{n}\end{pmatrix}=T_{n}\begin{pmatrix}s_{n-1}\\t_{n-1}\end{pmatrix},\begin{pmatrix}s_{0}\\ t_{0}\end{pmatrix}=\begin{pmatrix}1\\1\end{pmatrix}\end{equation}
where \begin{equation}\label{e.GZ2}T_{n}=\left\{\begin{aligned}\frac{1}{\rho_{n}}\begin{pmatrix}-\alpha_{n} &z^{-1}\\ z &-\overline{\alpha}_{n}\end{pmatrix} \qquad n \text{ is even,}\\\frac{1}{\rho_{n}}\begin{pmatrix}-\overline{\alpha}_{n} &1\\ 1&-\alpha_{n}\end{pmatrix}\qquad n\text{ is odd.}\end{aligned}\right.\end{equation}
Putting \eqref{e.szegorecnorm}, \eqref{e.GZ1}, and \eqref{e.GZ2} together for $z \in \partial \mathbb{D}$, a direct calculation gives
\begin{equation}\label{e.equinorm}
|s_{n}|=|\varphi_{n}|, \quad n\in\mathbb{Z}.
\end{equation}
Note that for $n > 0$, $\varphi_{-n}(z)$ is obtained via
\begin{equation}\label{e.Sz}
\begin{pmatrix} \varphi_{-n}(z) \\ \varphi^{*}_{-n}(z) \end{pmatrix} = \tilde{S}(\alpha_{-n+1},z)^{-1} \begin{pmatrix} \varphi_{-n+1}(z) \\ \varphi_{-n+1}^{*}(z) \end{pmatrix}, \; \begin{pmatrix} \varphi_{0}(z) \\ \varphi^{*}_{0}(z) \end{pmatrix} = \begin{pmatrix} 1 \\ 1 \end{pmatrix}.
\end{equation}
The $\varphi_{-n}(z)$'s are polynomials of $\frac{1}{z}$ with degree $n$ and $\varphi_{-n}^{*}(z)=z^{-n}\overline{\varphi_{n}(\frac{1}{\overline{z}})}$.
This relation inspires us to study the Szeg\H{o} cocycle instead of the Gesztesy-Zinchenko cocycle directly.

\subsection{Dynamically Defined Verblunsky Coefficients}

In this subsection we emphasize that the Verblunsky coefficients we study are obtained by sampling along the orbits of a discrete-time dynamical system, and hence they embed into the general theory of dynamically defined Verblunsky coefficients. We explain the primary objects of interest in this scenario and some of the key results that are known to hold. However, for the sake of simplicity we do not make the general framework explicit, but rather introduce all quantities for the specific case at hand. We merely point out that many of the quantities and results below make sense in a more general setting.

\subsubsection{The Density of States Measure}\label{CMV}

Averaging the spectral measure corresponding to the pair $(\mathcal{E}_x , \delta_0)$ with respect to the normalized Lebesgue measure on $\mathbb{T}^d$, we obtain the \textit{density of states measure} $dk$ on $\partial \mathbb{D}$:\footnote{We will sometimes view $\partial \mathbb{D}$ as $\R / (2 \pi \Z)$. In particular, when we talk about the weight assigned by the density of states measure to an arc on the unit circle, we typically denote that by $k(a,b)$, where $(a,b)$  is a real interval.}
$$
\int g \, dk = \int_{\mathbb{T}^d} \langle \delta_0, g(\mathcal{E}_x) \delta_0 \rangle \, dx
$$
for any $g\in C(\partial \D)$. Note that by shift-invariance of the Lebesgue measure on $\mathbb{T}^d$, $dk$ can alternatively be defined as one-half times the average of $\Lambda_x$.

Under suitable assumptions, the density of states measure can also be interpreted as the \textit{density of zeros measure}, which is the weak limit of the finitely supported probability measures obtained by placing point masses at the zeros of the orthogonal polynomials (according to their multiplicities).

Note that the zeros of the orthogonal polynomials lie in the open unit disk $\mathbb{D}$ and hence the existence of the density of zeros measure as a measure on the unit circle $\partial \mathbb{D}$ is a non-trivial property that will not always hold. For a very simple counterexample, note that the all zeros lie at the origin when the Verblunsky coefficients vanish, and hence the density of zeros measure exists in this case but is not a measure on the unit circle (rather it is the normalized point measure at the origin). On the other hand, the density of zeros exists as a measure on the unit circle when the function $\alpha$ defining the Verblunsky coefficients via \eqref{setting} satisfies
$$
\int_{\mathbb{T}^d} \ln (|\alpha(x)|) \, dx > - \infty,
$$
compare \cite[Theorem~10.5.19]{Simon2} and the general discussion in \cite[Section~10.5]{Simon2}. In view of \eqref{setting2}, this sufficient condition holds true for every $\lambda \in (0,1)$. In other words, we can use each of the two interpretations of the measure, as the density of states and as the density of zeros, in our setting.

\subsubsection{The Szeg\H{o} Cocycle, the Lyapunov Exponent, and the Rotation Number}

A \textit{quasi-periodic cocycle} $(\omega, A)\in \R^{d}\times C^0(\T^{d}, \mathrm{SL}(2,\C))$ is a linear skew product:
\begin{eqnarray*}
(\omega,A):&\T^{d} \times \C^2 \to \T^{d} \times \C^2\\
\nonumber &(x,v) \mapsto (x+\omega,A(x) \cdot v),
\end{eqnarray*}
where $\mathbb{T}^{d}=\mathbb{R}^{d}/\mathbb{Z}^{d}$ and $\omega\in\mathbb{T}^{d}$ is rationally independent. For $n \geq 1$, the products are defined as
$$
A^n(x)=A(x+(n-1)\omega) \cdots A(x),
$$
and $A^{-n}(x)=(A^n(x-n\omega))^{-1},$ then we can define \textit{Lyapunov exponent} as
$$
\gamma(\omega, A)=\lim_{n\rightarrow \infty} \frac {1} {n} \int \ln \|A^n(x)\| \, dx.
$$
Assume now that $A\in C^0(\T, \mathrm{SL}(2,\R))$ is homotopic to the identity. Then there exist $\psi:\T^{d} \times \T \to \R$ and $u:\T^{d} \times \T \to \R^+$ such that $$ A(x) \cdot \left (\bm \cos 2 \pi y \\ \sin 2 \pi y \em \right )=u(x,y) \left (\bm \cos 2 \pi (y+\psi(x,y)) \\ \sin 2 \pi (y+\psi(x,y)) \em \right ). $$ The function $\psi$ is called a {\it lift} of $A$.  Let $\mu$ be any probability measure on $\T^{d} \times \T$ which is invariant by the continuous map $T:(x,y) \mapsto (x+\omega,y+\psi(x,y))$, projecting over Lebesgue measure on the first coordinate (for instance, take $\mu$ as any accumulation point of $\frac {1} {n} \sum_{k=0}^{n-1} T_*^k \nu$, where $\nu$ is Lebesgue measure on $\T^{d} \times \T$). Then the number
$$
\rho(\omega,A) = \int \psi \, d\mu \mod \Z
$$
does not depend on the choices of $\psi$ and $\mu$, and is called the {\it fibered rotation number} of $(\omega,A)$, see \cite {H} and \cite {JM82}. For any $C\in \mathrm{SL}(2,\R),$ it is immediate from the definition that we have the following:
\begin{equation}\label{rotationnumberper}
|\rho(\omega, A)-\rho(\omega,C)|\leq\|A(x)-C\|^{\frac{1}{2}}_{C^{0}}.
\end{equation}

Once we have the definition of the fibered rotation number for $\mathrm{SL}(2,\mathbb{R})$ cocycles, we can parlay this concept to $\mathrm{SU}(1,1)$ cocycles, since
$$
M=\frac{1}{1+i}\begin{pmatrix}1 &-i\\ 1 &i\end{pmatrix}
$$
induces an isomorphism between $\mathrm{SL}(1,1)$ and $\mathrm{SL}(2,\mathbb{R})$:
$$
M^{-1} \mathrm{SU}(1,1) M = \mathrm{SL}(2,\mathbb{R}).
$$

The \textit{Szeg\H{o} cocycle} is defined by letting $\alpha_{n}(x)=\alpha(x+(n-1)\omega)$ and $A_{z}=S(\alpha,z)$ (with $S(\alpha,z)$ as in \eqref{e.renszegococyclemap}). Denote this cocycle by $(\omega, S(\alpha,z))$.

The rotation number $\rho(z)=\rho(\omega, S(\alpha,z))$ of the Szeg\H{o} cocycle and the DOS measure $dk$ associated with the family $\mathcal{E}_{x}$ are related by the following formula \cite[Theorem 8.3.3]{Simon2}:
\begin{equation}\label{rotationDOS}
2\rho(e^{i\zeta})=k(0,\zeta), \; \zeta\in [0,2\pi).
\end{equation}

We also denote the \textit{Lyapunov exponent} of the Szeg\H{o} cocycle and the renormalized Szeg\H{o} cocycle by
%\begin{equation}\label{LEnot}
$$\tilde{\gamma}(z)=\gamma(\omega,\tilde{S}(\alpha,z)), \; \gamma(z)=\gamma(\omega, S(\alpha,z)),$$
%\end{equation}
respectively. It follows readily from \eqref{e.renszegococyclemap} that
\begin{equation}\label{LErel}
\gamma(z) = \tilde{\gamma}(z)-\frac{1}{2}\log |z|.
\end{equation}
Note in particular that $\gamma(z)$ and $\tilde{\gamma}(z)$ coincide for $z \in \partial \mathbb{D}$.

The Lyapunov exponent and the DOS measure are related by the well known Thouless formula (cf.~\cite[Section~10.5]{Simon2}):
\begin{equation}\label{TH}
\tilde{\gamma}(z)=-\ln\rho_{\infty}+\int \ln|1-ze^{-i\theta}| \, dk(\theta),
\end{equation}
where
$$
\rho_{\infty} = \exp \left( \frac{1}{2} \int \log( 1-|\alpha(x)|^{2}) \, dx \right).
$$

\subsubsection{Uniform Hyperbolicity and the Set $\Sigma$}

We say that a cocycle $(\omega,A)$ is {\it uniformly hyperbolic} if there exist constants $c,\gamma>0$ and a continuous splitting $\mathbb{CP}^1=E^{s}(x)\oplus  E^{u}(x)$, $x \in \mathbb{T}^d$, such that
\begin{enumerate}

\item $A(x)E^{*}(x)=E^{*}(x+\omega)$ for $*=s,u$.

\item For $n\in \mathbb{Z}_{+},x\in X,\xi_{u}\in E^{u}(x),\xi_{s}\in E^{s}(x)$,
$$
\Vert A^{n}(x)\xi_{s}\Vert\leq ce^{-\gamma n}\Vert \xi_{s}\Vert,\qquad \Vert A^{-n}(x)\xi_{u}\Vert\leq ce^{-\gamma n}\Vert \xi_{u}\Vert.
$$
\end{enumerate}
We denote by $\mathcal{UH}$ the set of all uniformly hyperbolic cocycles.

The following result from \cite{DFLY} relates the spectrum of extended CMV matrices to their Szeg\H{o} cocycles:

\begin{Theorem}\label{uniformhyperbolic}
There is a compact set $\Sigma \subset \partial \mathbb{D}$ with $\sigma(\mathcal{E}_{x}) = \Sigma$ for every $x \in \mathbb{T}^d$. Moreover, this uniform spectrum $\Sigma$ is equal to $\partial \mathbb{D} \backslash U$, where
$$
U = \{ z \in \partial \mathbb{D} : (\omega,\tilde{S}(\alpha,z)) \in \mathcal{UH} \}.
$$
\end{Theorem}

We emphasize two important points. The spectrum in the two-sided case is independent of $x$ and it is described by an explicit dynamical property of the associated Szeg\H{o} cocycles. For this reason, it is beneficial to consider the two-sided case in the dynamically defined situation (assuming the underlying dynamical system is invertible) even if one is ultimately interested in the spectral analysis of the one-sided case. Namely, the two-sided spectrum $\Sigma$ is determined via the above characterization in terms of the absence of uniform hyperbolicity, and this set then in turn serves as the essential spectrum of the one-sided CMV matrix $\mathcal{C}_x$, which is $x$-independent as pointed out above. Moreover, the $x$-dependent discrete spectrum of $\mathcal{C}_x$ is quite tame: there is at most one point in each gap of $\Sigma$.

\section{Quantitative Almost Reducibility}\label{sec.3}

Our proof is based on quantitative almost reducibility results for $\mathrm{SU}(1,1)$ cocycles. We recall that $(\omega,A_1)$ is \textit{analytically conjugate} to $(\omega, A_2)$ if there exists $B \in C^{\omega}(2\T^{d}, \mathrm{SU}(1,1))$ such that
$$
B(\cdot+\omega) A_1(\cdot) B^{-1}(\cdot) = A_2(\cdot).
$$
The cocycle $(\omega, A)$ is said to be \textit{almost reducible} if the closure of its analytic conjugations contains a constant.

For any sufficiently small $\epsilon_{0} > 0$ and any $r > 0$, let us define the following sequences:
$$
\epsilon_{j} = \epsilon_{0}^{2^{j}}, \quad r_{j} = \frac{r}{2^{j}}, \qquad N_{j} = \frac{4^{j+1} \ln \epsilon_{0}^{-1}}{r}.
$$
Then we have the following quantitative almost reducibility results, which are based on the modified KAM scheme developed in \cite{CCYZ,LYZZ} (see Proposition~\ref{rd} in Appendix~\ref{sec.app} for the precise statement we employ).

\begin{Proposition}\label{rd1}
Assume that $\kappa, \tau, r > 0$ and $\omega \in \mathrm{DC}(\kappa,\tau)$. Let $S_{0} \in \mathrm{SU}(1,1)$, $f_{0} \in C^{\omega}_{r}(\mathbb{T}^{d}, \mathrm{su}(1,1))$ with
$$
\Vert f_{0} \Vert_{r} \leq \epsilon_{0} \leq \frac{D_{0}}{\left\| S_{0}\right\|^{C_{0}}} \left( \frac{r}{2} \right)^{C_{0}\tau},
$$
where $D_{0} = D_0(\kappa,\tau)$ and $C_{0}$ is a numerical constant. Then for any $j \geq 1$, there exists $B_{j} \in C^{\omega}_{r_{j}}(2\T^{d}, \mathrm{SU}(1,1))$ such that
$$
B_{j}(x+\omega)(S_{0}e^{f_{0}(x)}) B^{-1}_{j}(x) = S_{j} e^{f_{j}(x)},
$$
where $\Vert f_{j}(x) \Vert_{r_{j}}\leq \epsilon_{j}$ and $B_{j}$ satisfies
\begin{align}
\label{normOfB}
\left\| B_{j}\right\|_{0} & \leq \epsilon_{j-1}^{-\frac{1}{192}}, \\
\label{degreeOfB}\vert \deg{B_{j}}\vert & \leq 2N_{j-1}.
\end{align}

{\rm (a)}  For any $0 < \vert m \vert < N_{j-1}$, we denote
\begin{align}
\label{e.lambdamjdef}
\Lambda_{m}(j) = \{ z \in \Sigma : & \Vert 2 \rho( \omega, S_{j-1})- \langle m, \omega \rangle \Vert_{\mathbb{R}/\mathbb{Z}} < \epsilon_{j-1}^{\frac{1}{15}}  \}.
\end{align}
Then if $z \in \Lambda_{m}(j)$, we have the following precise expression:
$$
S_{j} = \exp \begin{pmatrix} i t_{j} & v_{j} \\ \bar{v}_{j} & -i t_{j} \end{pmatrix},
$$
where $t_{j} \in \mathbb{R}$, $v_{j} \in \mathbb{C}$, and $\vert t_{j} \vert \leq \epsilon_{j-1}^{\frac{1}{16}}$, $\vert v_{j} \vert \leq \epsilon_{j-1}^{\frac{15}{16}}$.

{\rm (b)}  Moreover, there always exist unitary matrices $U_{j}\in \mathrm{SL}(2,\mathbb{C})$ such that
\begin{equation}\label{e.conj}
U_{j} S_{j} e^{f_{j}(x)} U_{j}^{-1} = \begin{pmatrix} e^{2\pi i\rho_{j}} & c_{j} \\ 0 & e^{-2\pi i\rho_{j}} \end{pmatrix} + F_{j}(x)
\end{equation}
where $\rho_{j} \in \mathbb{R} \cup i \mathbb{R}$, with estimates $\Vert F_{j} \Vert_{r_{j}} \leq \epsilon_{j}$, and
\begin{equation}\label{Btimesc}
\Vert B_{j} \Vert_{0}^{2} \vert c_{j} \vert \leq 8 \Vert S_{0} \Vert.
\end{equation}
\end{Proposition}

\begin{pf}
We will prove  the result by iteration.

\textbf{First step:} By our choice of $\epsilon_{0}$, we can apply Proposition \ref{rd} once to $(\omega,S_{0} e^{f_{0}(x)})$ and obtain the following:

There exists $B_{1} \in C^{\omega}_{r_{1}}(2\T^{d},\mathrm{SU}(1,1))$ such that
$$
B_{1}(x+\omega) S_{0} e^{f_{0}(x)} B_{1}^{-1}(x) = S_{1} e^{f_{1}(x)}
$$
with the following estimates:
$$
\|f_{1}\|_{r_{1}} \leq \epsilon_{1}, \; \|B_{1}\|_{0} \leq \epsilon_{0}^{-\frac{1}{192}}, \; |\deg B_{1}| \leq 2N_{0}.
$$

\textbf{Inductive step:}
Suppose we have already completed the $j$-th step and are at the ($j+1$)-th step. That is,  there exists $B_{j} \in C^{\omega}_{r_{j}}(2\T^{d}, \mathrm{SU}(1,1))$ such that
$$
B_{j}(x+\omega)S_{0}e^{f_{0}(x)}B_{j}^{-1}(x)=S_{j}e^{f_{j}(x)}
$$
with estimates
$$
\|f_{j}\|_{r_{j}} \leq \epsilon_{j}, \; \|B_{j}\|_{0} \leq \epsilon_{j-1}^{-\frac{1}{192}}, \; |\deg B_{j}| \leq 2N_{j-1}.
$$
Now we consider the $(j+1)$-th step. By our choice of $\epsilon_{0}$, one can check that
\begin{equation}\label{ic}
\epsilon_{j}\leq \frac{D_{0}}{\left\|S_{j}\right\|^{C_{0}}}(r_{j}-r_{j+1})^{C_{0}\tau}.
\end{equation}
Indeed, $\epsilon_{j}$ on the left side of \eqref{ic} decays super-exponentially with $j$, while $(r_{j}-r_{j+1})^{C_{0}\tau}$ on the right side decays exponentially with $j$. Note that \eqref{ic} implies that Proposition~\ref{rd} can be applied iteratively, and hence we can construct
$$
\bar{B}_{j+1} \in C_{r_{j+1}}^{\omega}(2\mathbb{T}^{d}, \mathrm{SU}(1,1)), \; S_{j+1} \in \mathrm{SU}(1,1), \; f_{j+1} \in C_{r_{j+1}}^{\omega}(\mathbb{T}^{d}, \mathrm{su}(1,1))
$$
such that
$$
\bar{B}_{j+1}(x+\omega)S_{j}e^{f_{j}(x)}\bar{B}^{-1}_{j+1}(x)=S_{j+1}e^{f_{j+1}(x)}.
$$
Let $B_{j+1}=\bar{B}_{j+1}B_{j}$, then we have $$B_{j+1}(x+\omega)S_{0}e^{f_{0}(x)}B^{-1}_{j+1}(x)=S_{j+1}e^{f_{j+1}(x)}.$$
In order to verify \eqref{normOfB} and \eqref{degreeOfB} for $B_{j+1}$, we need to distinguish the following two cases:

\textbf{Non-resonant case:} If for any $n\in \mathbb{Z}^{d}$ with $0<\vert n\vert\leq N_{j}$ we have
$$
\left\| 2\rho(\omega, S_j)-\langle n,\omega \rangle\right\|_{\mathbb{R}/ \mathbb{Z}}\geq \epsilon_{j}^{\frac{1}{15}},
$$
then by the non-resonant case of Proposition \ref{rd},
$$
\left\|\bar{B}_{j+1}-id\right\|_{r_{j+1}} \leq \epsilon_{j}^{\frac{1}{2}}, \; \left\|f_{j+1}\right\|_{r_{j+1}} \leq \epsilon_{j+1}.
$$
In this case we have $\deg{B_{j+1}}=\deg{B_{j}}$, since $\bar{B}_{j+1}$ is homotopic to the identity. Therefore,
\begin{align*}
%\left\| B_{j+1}\right\|_{r_{j+1}} & \leq (1+\epsilon_{j}^{\frac{1}{2}})\Vert B_{j}\Vert_{r_{j}}, \\
\Vert B_{j+1} \Vert_{0} & \leq (1+\epsilon_{j}^{\frac{1}{2}})\epsilon_{j-1}^{-\frac{1}{192}}\leq \epsilon_{j}^{-\frac{1}{192}}, \\
\vert \deg{B_{j+1}} \vert & = \vert\deg{B_{j}} \vert \leq 2N_{j-1}\leq 2N_{j}.
%\Vert S_{j+1}-S_{j} \Vert & \leq 2\epsilon_{j}.
\end{align*}
These verify \eqref{normOfB} and \eqref{degreeOfB} for $B_{j+1}$ in this case.
\par
\textbf{Resonant case:} If $z \in \Lambda_{m}(j+1)$ for some $m \in \mathbb{Z}^{d}$ with $0 < \vert m \vert \leq N_{j}$, then by the resonant case of Proposition \ref{rd}, we have $\deg \bar{B}_{j+1}=m$ and
\begin{align*}
\nonumber \Vert \bar{B}_{j+1}\Vert_{0} & \leq C\vert m\vert^{\tau}\leq\epsilon_{j}^{-\frac{1}{384}}, \\
\nonumber \left\| f_{j+1}\right\|_{r_{j+1}} & \leq \epsilon_{j} e^{-r_{j+1} \epsilon_{j}^{-\frac{1}{18\tau}}} \leq \epsilon_{j+1}.
\end{align*}
Moreover, we can write $S_{j+1}=\exp \begin{pmatrix}it_{j+1} &v_{j+1}\\\bar{v}_{j+1} &-it_{j+1}\end{pmatrix}$, where $t_{j+1}\in\mathbb{R},v_{j+1}\in\mathbb{C}$ with $|t_{j+1}|\leq\epsilon_{j}^{\frac{1}{16}},|v_{j+1}|\leq\epsilon_{j}^{\frac{15}{16}}.$
Therefore,
\begin{align*}
\|B_{j+1}\|_{0} & \leq \epsilon_{j-1}^{-\frac{1}{192}}\epsilon_{j}^{-\frac{1}{384}}\leq\epsilon_{j}^{-\frac{1}{192}}, \\
\vert \deg B_{j+1} \vert & \leq 2N_{j-1}+N_{j}\leq 2N_{j}.
\end{align*}

By induction we have proven \eqref{normOfB} and \eqref{degreeOfB} for each $j\geq 1$. As a consequence of the resonant case, the statement of $(a)$ was verified. We are left to prove the statement $(b)$.

From the above iteration, we know after each resonant step $j$, one can write $S_{j}=\exp S_{j}''$ with $\|S_{j}''\|\leq\epsilon_{j-1}^{\frac{1}{16}}$. Let $U_{j}\in \mathrm{SL}(2,\mathbb{C})$ such that \eqref{e.conj} holds, then
\begin{equation}\label{e.estc}|c_{j}|\leq 2\|S_{j}''\|\leq 2\epsilon_{j-1}^{\frac{1}{16}}.\end{equation}

Suppose there are two resonance sites $n_{j_{i}},n_{j_{i+1}}$, which happen at the KAM steps $j_{i}+1,j_{i+1}+1$, respectively. Then we want to show that $\vert n_{j_{i+1}} \vert \geq \epsilon_{j_{i}}^{-\frac{1}{18\tau}} \vert n_{j} \vert$. By the resonance condition of the $(j_{i+1}+1)$-th step, $\vert \rho_{j_{i+1}}-\frac{\langle n_{j_{i+1}},\omega\rangle}{2}\vert\leq\frac{1}{2}\epsilon^{\frac{1}{15}}_{j_{i+1}}$, hence $\vert\rho_{j_{i+1}}\vert>\frac{\kappa}{3\vert n_{j_{i+1}}\vert^{\tau}}$. On the other hand,  according to Proposition~\ref{rd}, after the $(j_{i}+1)$-th step, $\vert \rho_{j_{i}+1} \vert \leq \epsilon_{j_{i}}^{\frac{1}{16}}$. Thus we have
$$
\vert n_{j_{i+1}} \vert \geq \epsilon_{j_{i}}^{-\frac{1}{18\tau}}\vert n_{j_{i}}\vert.
$$
This implies that the resonant steps are actually very far from each other.

Since the steps between $j_{i}+1$ and $j_{i+1}+1$ are all non-resonant, according to the non-resonant case of Proposition \ref{rd}, we have the following:
For $j_{i}+1 < j \leq j_{i+1}$,
\begin{equation}\label{e.nondiff}
\left\| S_{j} - S_{j-1} \right\| \leq 2\epsilon_{j-1}.
\end{equation}
Since the $(j_{i}+1)$-th step is resonant, we have
$$
S_{j_{i}+1} =\exp{\begin{pmatrix} it_{j_{i}+1} &v_{j_{i+1}}\\ \bar{v}_{j_{i}+1} &-it_{j_{i}+1} \end{pmatrix}}
$$
with $\vert t_{j_{i}+1} \vert \leq \epsilon_{j_{i}}^{\frac{1}{16}}$ and $\vert v_{j_{i}+1} \vert \leq \epsilon_{j_{i}}^{\frac{15}{16}}$. There exists a unitary matrix $U_{j_{i}+1} \in \mathrm{SL}(2,\mathbb{C})$ such that
$$
U_{j_{i}+1} S_{j_{i}+1} e^{f_{j_{i}+1(x)}} U_{j_{i}+1}^{-1} = \begin{pmatrix} e^{ 2\pi i\rho_{j_{i}+1}} & c_{j_{i}+1} \\ 0 & e^{- 2\pi i\rho_{j_{i}+1}} \end{pmatrix} + F_{j_{i}+1}(x),
$$
where $\vert c_{j_{i}+1} \vert\leq 2\|S_{j_{i}+1}''\| \leq 2\epsilon_{j_{i}}^{\frac{1}{16}}$ and $\|F_{j_{i}+1}\|_{r_{j_{i}+1}}\leq\epsilon_{j_{i}+1}$. Let $S_{j} = \exp S^{''}_{j}$. By \eqref{e.nondiff}, we have the following:
$$
\Vert S^{''}_{j} \Vert \leq 2\epsilon_{j_{i}}^{\frac{1}{16}} + \sum^{j}_{k=j_{i} + 2} \Vert S_{k} - S_{k-1} \Vert \leq 4\epsilon_{j_{i}}^{\frac{1}{16}}.
$$
In this case, we have
\begin{equation}\label{nonresonantc}
\vert c_{j} \vert \leq 2\Vert S_{j}^{''} \Vert \leq 8\epsilon_{j_{i}}^{\frac{1}{16}}.
\end{equation}
Since $\Vert B_{j_{i}+1} \Vert_{0} \leq \epsilon_{j_{i}}^{-\frac{1}{192}}$ and all steps between $j_{i}+1$ and $j$ (including $j$) are non-resonant, we have
\begin{equation}\label{nonresonantnorm}
\Vert B_{j} \Vert_{0} \leq \Pi_{k=j_{i}+1}^{j}(1+\epsilon_{k}^{\frac{1}{2}}) \epsilon_{j_{i}}^{-\frac{1}{192}} \leq 2\epsilon_{j_{i}}^{-\frac{1}{192}}.
\end{equation}

From the above argument and \eqref{normOfB},\eqref{e.estc}, we can see that after each resonant step $j$, we have $$
\Vert B_{j} \Vert_{0}^{2} \vert c_{j} \vert \leq 2\epsilon_{j-1}^{\frac{5}{96}}.
$$
While after each non-resonant step $j$, we are able to utilize the estimates of the last resonant step $j_{i}+1$ if it exists. By \eqref{nonresonantc} and \eqref{nonresonantnorm} we have
$$
\Vert B_{j} \Vert_{0}^{2} \vert c_{j} \vert \leq 32\epsilon_{j_{i}}^{\frac{5}{96}}.
$$
However, it is possible that no resonant steps happened within the first $j$ steps. In this case, each step is non-resonant and thus we can use the estimate $\|\bar{B_{i}}\|_{0}\leq 1+\epsilon_{i-1}^{\frac{1}{2}}$ for each $i\leq j$ and obtain
$$
\|B_{j}\|_{0} \leq \Pi_{i\geq 0}(1+\epsilon_{i}^{\frac{1}{2}}) \leq 2.
$$
Since $S_{j}-S_{0}=\sum_{i=1}^{j}(S_{i}-S_{i-1})$ and \eqref{e.nondiff}, we have $\|S_{j}\| \leq \|S_{0}\| + 2\epsilon_{0} \leq 2\|S_{0}\|$. This finishes the proof.
\end{pf}

Let
\begin{equation}\label{e.Kjdef}
K_{j}=\bigcup_{0<\vert m\vert\leq N_{j-1}}\Lambda_{m}(j)
\end{equation}
with $\Lambda_{m}(j)$ from \eqref{e.lambdamjdef}. Then for any $z\in K_{j}$, we have the following estimate for the growth of the transfer matrix $A_{z}^{s}$:

\begin{Corollary}\label{growth01}
For $z \in K_{j}$,  we have
$$
\sup \limits_{0 \leq s \leq C\epsilon_{j-1}^{-\frac{1}{16}}} \Vert A_{z}^{s} \Vert_{0} \leq C \epsilon_{j-1}^{-\frac{1}{96}}.
$$
\end{Corollary}

\begin{pf}
Since $z \in K_{j}$, by Proposition~\ref{rd1} there exists $B_{j}$ such that
$$
B_{j}(x+\omega) S_{0} e^{f_{0}(x)} B^{-1}_{j}(x) = S_{j} e^{f_{j}(x)},
$$
where
$$
S_{j} = \exp \begin{pmatrix} i t_{j} & \bar{v}_{j} \\ v_{j} & -i t_{j} \end{pmatrix}
$$
and
$$
\vert t_{j} \vert \leq \epsilon_{j-1}^{\frac{1}{16}}, \; \vert v_{j} \vert \leq \epsilon_{j-1}^{\frac{15}{16}}, \; \Vert f_{j} \Vert_{0} \leq \epsilon_{j}.
$$
Moreover, we have $\Vert S_{j} \Vert \leq 1 + \epsilon_{j-1}^{\frac{1}{16}}$. Since we have $\Vert B_{j} \Vert_{0} \leq \epsilon_{j-1}^{-\frac{1}{192}}$ by \eqref{normOfB}, this implies that
$$
\sup \limits_{0 \leq s \leq C \epsilon_{j-1}^{-\frac{1}{16}}} \Vert A_{z}^{s}\Vert_{0}\leq C\|B_{j}\|^{2}_{0}\leq C\epsilon_{j-1}^{-\frac{1}{96}},
$$
concluding the proof.
\end{pf}

We also have the following:

\begin{Lemma}\label{LM4}
For $z\in K_{j}$, there exists  $n_{j}\in\mathbb{Z}^{d}$ with $\vert n_{j}\vert \leq 2N_{j-1}$ such that
$$
\|2\rho(\omega,S_{0}e^{f_{0}(x)})-\langle n_{j},\omega\rangle\|_{\mathbb{R}/\mathbb{Z}}\leq 2\epsilon_{j-1}^{\frac{1}{15}}.
$$
\end{Lemma}

\begin{pf}
For $z\in K_{j}$, we have
$$
\Vert 2\rho_{j-1}-\langle n^{*},\omega\rangle\Vert_{\mathbb{R}/\mathbb{Z}}\leq\epsilon_{j-1}^{\frac{1}{15}}
$$
for some $n^{*}\in\mathbb{Z}^{d}$ with $0 < \vert n^{*} \vert \leq N_{j-1}$. Since $B_{j-1}(\cdot)$ conjugates $(\omega,S_{0}e^{f_{0}(\cdot)})$ into $(\omega,S_{j-1}e^{f_{j-1}(\cdot)})$, we have
$$
2\rho(\omega,S_{0}e^{f_{0}(x)})+\langle \deg B_{j-1},\omega\rangle=2\rho(\omega,S_{j-1}e^{f_{j-1}(x)}).
$$
By \eqref{rotationnumberper} and $\|f_{j-1}\|_{0}\leq\epsilon_{j-1}$, we have
$$
\vert \rho(\omega,S_{j-1}e^{f_{j-1}(x)})-\rho_{j-1}\vert\leq c\epsilon_{j-1}^{\frac{1}{2}}.
$$
These together give
$$
| 2\rho(\omega,S_{0}e^{f_{0}(x)})+\langle\deg B_{j-1},\omega\rangle-\langle n^{*},\omega\rangle|\leq2\epsilon_{j-1}^{\frac{1}{15}}.
$$
Let $n_{j}=-\deg B_{j-1}+n^{*}$. Then
$$
\| 2\rho(\omega,S_{0}e^{f_{0}(x)})-\langle n_{j},\omega\rangle\|_{\mathbb{R}/ \mathbb{Z}}\leq 2\epsilon_{j-1}^{\frac{1}{15}}.
$$
Moreover, by \eqref{degreeOfB} of Proposition $\ref{rd1},$ we have $|\deg B_{j-1}|\leq 2N_{j-2}$, and this implies that $|n_{j}|\leq 2N_{j-2}+N_{j-1}\leq 2N_{j-1}$, as claimed.
\end{pf}

\section{A Lower Bound for the DOS}\label{sec.4}

In this section we consider CMV matrices with analytic quasi-periodic Verblunsky coefficients of the form
$$
\alpha_{n}(x) = \alpha(x+(n-1)\omega), \; n \in \mathbb{Z}, \; x \in \T^{d}, \; \omega \in \mathrm{DC}(\kappa,\tau),
$$
where
$$
\alpha(x) = \lambda e^{2\pi i h(x)}, \; h(x)\in C^{\omega}(\mathbb{T}^{d},\mathbb{R}), \; \lambda\in(0,1).
$$
Since the $x$-independent spectrum $\Sigma$ of $\mathcal{E}_x$ is a subset of $\partial\mathbb{D}$, we will use the notation $e^{i\zeta} \in \Sigma$ to make this explicit.

We will show that if $\lambda$ is small enough, then we can establish (see Lemma~\ref{LB} below) a lower bound for the DOS measure for $\mathcal{C}_{x}$ of the form
\begin{equation}\label{lower}
k(\zeta -\epsilon,\zeta+\epsilon)\geq c\epsilon^{\frac{3}{2}}, e^{i\zeta}\in\Sigma,
\end{equation}
which will be a key ingredient in the subsequent considerations.

If $|\lambda|$ is small enough, then $S(\alpha,z)$ is close to a constant, and hence employing Proposition~\ref{rd1} we will be able to show (see Corollary~\ref{upper} below) that the DOS measure is $\frac{1}{2}$-H\"older continuous on $\Sigma$. Once this has been accomplished, \eqref{lower} follows from the well-known Thouless formula.

Denoting $\rho = \sqrt{1-\vert \alpha\vert^{2}}=\sqrt{1-\lambda^{2}}$, we have the following decomposition
$$
S(\alpha,z)=\begin{pmatrix} z^{\frac{1}{2}} & 0 \\ 0 & z^{-\frac{1}{2}} \end{pmatrix} \rho^{-1} \begin{pmatrix} 1 & -\overline{\alpha}z^{-1} \\ -\alpha z & 1 \end{pmatrix}.
$$
Let
$$
S_{0}(\alpha,z) = \begin{pmatrix} z^{\frac{1}{2}}  & 0 \\ 0 & z^{-\frac{1}{2}} \end{pmatrix}
$$
and
$$
P(x) = \rho^{-1}\begin{pmatrix} 1 & -\lambda e^{-2\pi ih(x)}z^{-1} \\ -\lambda e^{2\pi ih(x)}z & 1 \end{pmatrix}.
$$
Since
$$
P(x) - id = \frac{1}{\rho} \begin{pmatrix} 1 - \rho & -\lambda e^{2\pi ih(x)} z^{-1} \\ -\lambda e^{2\pi i h(x)} z & 1 - \rho \end{pmatrix}
$$
for $z \in \partial \mathbb{D}$, we have
$$
\|P(x) - id\|_{r} \leq 2|\lambda| e^{2\pi |h|_{r}}.
$$
Since the tangent space $T_{id} \mathrm{SU}(1,1)$ is the algebra $\mathrm{su}(1,1)$, the exponential map sets up a local diffeomorphism between the neighborhood of $id \in \mathrm{SU}(1,1)$ and that of $0 \in \mathrm{su}(1,1)$. This means that if $\lambda$ is sufficiently small, we have $f_{0}(x) \in C_{r}^{\omega}(\mathbb{T}^{d},\mathrm{su}(1,1))$ such that $e^{f_{0}(x)}=P(x)$, and hence we have $S(\alpha,z)=S_{0}(\alpha ,z)e^{f_{0}(x)}$. The smallness of $\lambda$ can be viewed as the smallness of $\Vert f_{0}\Vert_{r}$. If we further assume that $\lambda$ satisfies the condition
\begin{equation}\label{e.lambdasmall}
|\lambda| \leq \frac{D_{0}}{2e^{2\pi |h|_{r}}} \left( \frac{r}{2} \right)^{C_{0}\tau},
\end{equation}
where the constants $C_{0}, D_{0}$ are consistent with Proposition~\ref{rd}, then we can apply Proposition~\ref{rd1} to the system $(\omega, S_{0}e^{f_{0}(x)})$, and we obtain the following consequences.

\begin{Lemma}\label{Prop3.2}
If $\lambda$ obeys \eqref{e.lambdasmall}, then for any $e^{i\zeta} \in \Sigma$ and $\epsilon>0$, we have
\begin{equation}\label{zero}
\gamma(e^{i\zeta}) = 0
\end{equation}
and
\begin{equation}\label{growth}
\gamma((1+\epsilon) e^{i\zeta}) \leq C \epsilon^{\frac{1}{2}}.
\end{equation}
\end{Lemma}

\begin{pf}
To prove \eqref{zero}, it suffices to show that if $e^{i\zeta} \in \Sigma$, the cocycle grows at most linearly, that is,
\begin{equation}\label{lineargrow}
\Vert A_{e^{i\zeta}}^{n}\Vert_{0} \leq Cn, \quad n \ge 1.
\end{equation}
We can apply Proposition \ref{rd1} and distinguish two cases:\\

\noindent  \textbf{Case 1:}
If $(\omega,S_{0}e^{f_{0}(x)})$ is reducible, then there exists $B\in C^{\omega}(2\mathbb{T}^{d},\mathrm{SU}(1,1))$ such that
$$
B(x+\omega)S_{0}e^{f_{0}(x)}B^{-1}(x)=S.
$$
Since $e^{i\zeta}\in \Sigma$, by Theorem \ref{uniformhyperbolic}, the cocycle $(\omega,S_{0}e^{f_{0}(x)})$ is not uniformly hyperbolic, and hence we have
$$
S = \begin{pmatrix} e^{2\pi i\rho} & c\\ 0 &e^{-2\pi i\rho}\end{pmatrix}
$$
with $\rho \in \mathbb{R}$. This implies that $\Vert A_{e^{i\zeta}}^{n}\Vert_{0} \leq \Vert B \Vert_{0}^{2}(1 + cn) \leq Cn$ for $n \ge 1$. \\

\noindent  \textbf{Case 2:}
If $(\omega,S_{0}e^{f_{0}(x)})$ is not reducible but almost reducible,
we need the following claim in order to describe the growth of the cocycle more precisely:
\begin{Claim}\label{hyperbolic}
Suppose $e^{i\zeta}\in\Sigma$, then for each $j>0$, there exists $\tilde{B}_{j}(x)$ such that
\begin{equation}\label{eq12062}
\tilde{B}_{j}(x+\omega)(S_{0}e^{f_{0}(x)})\tilde{B}^{-1}_{j}(x)=\begin{pmatrix} e^{2\pi i\rho_{j}} &c_{j}\\ 0 & e^{-2\pi i\rho_{j}}  \end{pmatrix}+\tilde{F}_{j}(x),
\end{equation}
with $\rho_{j} \in \mathbb{R},\|\tilde{B}_{j}\|_{0}\leq \epsilon_{j-1}^{-\frac{1}{192}}, \|\tilde{F}_{j}\|_{0} \leq \epsilon_{j}^{\frac{1}{4}}, \|\tilde{B}_{j}\|_{0}^{2}|c_{j}|\leq 8\|S_{0}\|$.
\end{Claim}

\begin{pf}
By Proposition \ref{rd1}, we obtain $B_{j}(x)$ such that
$$
B_{j}(x+\omega) S_{0} e^{f_{0}(x)} B^{-1}_{j}(x) = \begin{pmatrix} e^{2\pi i\rho_{j}} & c_{j} \\ 0 & e^{-2\pi i\rho_{j}} \end{pmatrix} + F_{j}(x),
$$
with $\rho_{j}\in\mathbb{R}\cup i\mathbb{R},\| B_{j}\|_{0}\leq\epsilon_{j-1}^{-\frac{1}{192}},\|F_{j}\|_{0}\leq \epsilon_{j}$ and $\|B_{j}\|_{0}^{2}|c_{j}|\leq 8\|S_{0}\|$. If $\rho_{j} \in \mathbb{R}$, then let $\tilde{B}_{j}(x) = B_{j}(x)$ and $\tilde{F}_{j}(x) = F_{j}(x)$, the claim follows immediately. Assume that $\rho_{j} \in i \mathbb{R}$ with $\vert \rho_{j}\vert> \epsilon_{j}^{\frac{1}{4}}$ and let
$$
Q_{j} = \begin{pmatrix} q_{j} & 0 \\ 0 & q_{j}^{-1} \end{pmatrix},
$$
where $q_{j} = \Vert B_{j}\Vert_{0}\epsilon_{j}^{\frac{1}{4}}$. Then we have
$$
Q_{j} \left[ \begin{pmatrix}e^{2\pi i\rho_{j}}&c_{j}\\0&e^{-2\pi i\rho_{j}}\end{pmatrix}+F_{j}(x) \right] Q_{j}^{-1}=\begin{pmatrix}e^{2\pi i\rho_{j}} & 0 \\ 0 & e^{-2\pi i \rho_{j}}\end{pmatrix}+F'_{j}(x),
$$
where
$$
F'_{j}(x) = Q F_{j} Q^{-1} + \begin{pmatrix} 0 & c_{j} q_{j}^{2} \\ 0 & 0 \end{pmatrix}.
$$
Since $|c_{j}|q_{j}^{2}=|c_{j}|\|B_{j}\|_{0}^{2}\epsilon_{j}^{\frac{1}{2}}$ and $\|Q_{j}F_{j}(x)Q_{j}^{-1}\|_{0}\leq \epsilon_{j}^{\frac{1}{2}}$,
we have  $\|F'_{j}\|_{0}\leq C\epsilon_{j}^{\frac{1}{2}}$.
We want to show that this implies the system is uniformly hyperbolic. Given a non-zero vector $(a,b)^T \in\mathbb{R}^{2}$ with $|a|\geq|b|$, let
$$
\begin{pmatrix} a' \\ b' \end{pmatrix} = \left[\begin{pmatrix} e^{2\pi i\rho_{j}} & 0 \\ 0 & e^{-2\pi i \rho_{j}} \end{pmatrix} + F'_{j}(x) \right] \begin{pmatrix} a \\ b \end{pmatrix} = \begin{pmatrix} e^{2\pi i \rho_{j}} a \\ e^{-2\pi i\rho_{j}} b \end{pmatrix} + F'_{j}(x) \begin{pmatrix}a\\ b\end{pmatrix}.
$$
Without loss of generality, assume $ 2\pi i\rho_{j}>0$. This implies
\begin{align*}
|a'| & \geq (e^{2\pi i\rho_{j}}-2\|F'_{j}(x)\|_{0})|a| \\
|b'| & \leq e^{-2\pi i\rho_{j}}|b|+2\|F'_{j}(x)\|_{0}|a|.
\end{align*}
Therefore, $|a'|-|b'|\geq (4\pi i\rho_{j}-2C\epsilon_{j}^{\frac{1}{2}})|a|\geq 0$. By the cone field criterion (compare, e.g., \cite{Yoc}), this implies the uniform hyperbolicity of $(\omega,S(\alpha,e^{i\zeta}))$, which conflicts with our assumption that $e^{i\zeta} \in \Sigma$, again by Theorem~\ref{uniformhyperbolic}.
So we have $\vert \rho_{j}\vert\leq\epsilon_{j}^{\frac{1}{4}}$ and we put it into the perturbation to obtain the following:
$$
\tilde{B}_{j}(x+\omega)(S_{0}e^{f_{0}(x)})\tilde{B}^{-1}_{j}(x)=\begin{pmatrix} 1 & c_{j}\\ 0 & 1  \end{pmatrix}+\tilde{F}_{j}(x),\|\tilde{F}_{j}\|_{0}\leq\epsilon_{j}^{\frac{1}{4}},
$$
where $\tilde{B}_{j}(x)=B_{j}(x)$.
This reduces to the case $\rho_{j}=0$ and Claim~1 is proved, as we indeed have $\rho_{j} \in \mathbb{R}$.
\end{pf}

In order to control the growth of the cocycle, we need the following lemma proved by Avila-Fayad-Krikorian \cite{AFK}:

\begin{Lemma}\label{afk}
We have that
$$
M_l(\id+\xi_l)\cdots M_0(\id+\xi_0)=M^{(l)}(\id+\xi^{(l)}),
$$
where $M^{(l)}=M_l\cdots M_0$ and
$$
\lVert \xi^{(l)}\rVert\leqslant e^{\sum_{k=0}^{l}\lVert M^{(k)}\rVert^2\lVert \xi_k\rVert}-1.
$$
\end{Lemma}
By the result of Claim \ref{hyperbolic}, let $M_k=S_{j}=\begin{pmatrix}e^{2\pi i\rho_{j}} &c_{j}\\ 0 &e^{-2\pi i\rho_{j}}\end{pmatrix}$, $\xi_k = S_{j}^{-1} \tilde{F}_{j}(x+k\omega)$, and apply Lemma~\ref{afk} to obtain
$$
A_{e^{i\zeta}}^{n}=\tilde{B}_{j}(x+n\omega)S^{n}_{j}(id+\xi^{(n)})\tilde{B}^{-1}_{j}(x),
$$
where $\|\xi^{(n)}\|\leq e^{\sum_{k=1}^{n}\|S^{k}_{j}\|^{2}\|\tilde{F}_{j}\|_{0}}-1$. Since $\rho_{j}\in\mathbb{R}$, we have $\|S_{j}^{k}\|\leq 1+k|c_{j}|$. These together with $\|\tilde{F}_{j}\|_{0}\leq \epsilon_{j}^{\frac{1}{4}}$ give
$$
\Vert A_{e^{i\zeta}}^{n}\Vert_{0}\leq \Vert \tilde{B}_{j}\Vert^{2}_{0}(1+n\vert c_{j}\vert)e^{\sum_{k=1}^{n}(1+k|c_{j}|)^{2}\epsilon_{j}^{\frac{1}{4}}}\leq \Vert \tilde{B}_{j}\Vert^{2}_{0}(1+n\vert c_{j}\vert)e^{n^{3}\epsilon_{j}^{\frac{1}{4}}}.
$$
Therefore, by \eqref{Btimesc} and the fact that $\|S_{0}\|\leq 1$ we have the following:
\begin{equation}\label{eq12063}
\sup_{0<n<\epsilon_{j-1}^{-\frac{1}{16}}}\Vert A_{e^{i\zeta}}^{n}\Vert_{0}\leq 2\Vert \tilde{B}_{j}\Vert^{2}_{0}(1+n\vert c_{j}\vert)\leq 2\epsilon_{j-1}^{-\frac{1}{96}}+16n.
\end{equation}
Since $\epsilon_{j}=\epsilon_{0}^{2^{j}}$, for any fixed large $n\in\mathbb{Z}_{+}$, there exists $j$ such that $n\in[\epsilon_{j-1}^{-\frac{1}{96}},\epsilon_{j-1}^{-\frac{1}{16}}]$, then by \eqref{eq12063}, we have $\Vert A_{e^{i\zeta}}^{n}\Vert_{0}\leq Cn,$ and this verifies \eqref{lineargrow}.

Since
$$
\gamma (e^{i\zeta}) = \lim\limits_{n \to \infty} \frac{1}{n} \int_{x \in \mathbb{T}^{d}} \ln \Vert A_{e^{i\zeta}}^{n}(x) \Vert \, dx,
$$
we obtain \eqref{zero} for any $e^{i\zeta}\in\Sigma.$

Let $I_{j} = (\epsilon_{j}^{\frac{1}{4}},\epsilon_{j-1}^{\frac{1}{48}})$. Then, for any small $\epsilon>0$, there exists $j$ such that $\epsilon\in I_{j}$. To prove \eqref{growth}, we need the following result:
\begin{Claim}\label{claim2}
There exists $W : 2\T^{d} \to \mathrm{SU}(1,1)$ analytic with $\left\| W\right\|_{0}\leq C\epsilon^{-\frac{1}{4}}$ such that
$$
Q(x) = W(x+\omega)S_{0}e^{f_{0}(x)}W^{-1}(x)
$$
satisfies the estimate
\begin{equation}\label{qqqqq}
\left\| Q\right\|_{0}\leq 1+C\epsilon^{\frac{1}{2}}.
\end{equation}
\end{Claim}

\begin{pf}
Let $\tilde{B}_{j}(x), \tilde{F}_{j}(x)$ be as above in Claim \ref{hyperbolic}, and let
$$
D = \begin{pmatrix} d & 0 \\ 0 & d^{-1} \end{pmatrix},
$$
where $d = \left\| \tilde{B}_{j}\right\|_{0}\epsilon^{\frac{1}{4}}$,
let $W(x)=D\tilde{B}_{j}(x)$. Since $\epsilon\in I_{j}$, we have $d\leq 1$, which implies $\left\| W\right\|_{0}\leq C\epsilon^{-\frac{1}{4}}$. By \eqref{eq12062} we have
$$
W(x+\omega)S_{0}e^{f_{0}(x)}W^{-1}(x) = \begin{pmatrix} e^{2\pi i\rho_{j}} & 0 \\ 0 & e^{-2\pi i\rho_{j}} \end{pmatrix}+\tilde{F}'_{j}(x),
$$
where
$$
\tilde{F}'_{j}(x)=\begin{pmatrix}0 &\left\| \tilde{B}_{j}\right\|^{2}_{0}\epsilon^{\frac{1}{2}}c_{j} \\ 0 & 0 \end{pmatrix} + D\tilde{F}_{j}(x)D^{-1},
$$
since $\|\tilde{B}_{j}\|_{0}^{2}|c_{j}|\leq 8\|S_{0}\|$ according to \eqref{Btimesc}, we have $\| \tilde{F}'_{j}\|_{0}\leq C\epsilon^{\frac{1}{2}}$ and $\rho_{j}\in \mathbb{R}$. Then, with $Q(x)=W(x+\omega)(S_{0}e^{f_{0}(x)})W^{-1}(x)$, \eqref{qqqqq} follows immediately. This completes the proof of Claim~2.
\end{pf}

Let $A=S(\alpha,e^{i\zeta}),B=S(\alpha,(1+\epsilon)e^{i\zeta})$, then $\Vert A-B\Vert\leq \epsilon.$ By the above argument, there exists $W: 2\T^{d}\to \mathrm{SU}(1,1)$ such that
$$
C=W(x+\omega)BW^{-1}(x)=W(x+\omega)AW^{-1}(x)+W(x+\omega)(B-A)W^{-1}(x).
$$
Thus we have $\|C\|_{0}\leq \| Q\|_{0}+C\epsilon^{\frac{1}{2}}$, by Claim \ref{claim2} and our choice of $W$,
$\Vert C\Vert_{0}\leq 1+C\epsilon^{\frac{1}{2}}.$ It follows that $\gamma((1+\epsilon)e^{i\zeta})\leq \ln \| C \|_{0}\leq C\epsilon^{\frac{1}{2}}$. This concludes the proof of Lemma~\ref{Prop3.2}.
\end{pf}

As a direct corollary, we have the following:

\begin{Corollary}\label{upper}
The DOS measure is $\frac{1}{2}$-H\"older continuous, that is for any $e^{i\zeta}\in \partial\mathbb{D}$ and $\epsilon>0$ small enough, we have
\begin{equation}\label{dos}
k(\zeta-\epsilon,\zeta+\epsilon)\leq C\epsilon^{\frac{1}{2}}.
\end{equation}
\end{Corollary}

\begin{pf}
Since the DOS measure is supported by $\Sigma$, we can limit our attention to the case $e^{i\zeta}\in \Sigma$.
By the Thouless formula \eqref{TH}, we have
$$
\tilde{\gamma}(re^{i\zeta})=-\ln\rho_{\infty} + \int \ln \vert 1-re^{i(\zeta-\theta)}\vert \, dk(\theta).
$$
For any $\epsilon$ which is small enough, using the fact $\tilde{\gamma}(e^{i\zeta}) = \gamma(e^{i\zeta}) = 0$ by Lemma $\ref{Prop3.2}$, we have
$$\begin{aligned}
\tilde{\gamma}((1+\epsilon)e^{i\zeta}) &=\tilde{\gamma}((1+\epsilon)e^{i\zeta})-\tilde{\gamma}(e^{i\zeta})\\
&=\int \ln \left| \frac{1-(1+\epsilon)e^{i(\zeta-\theta)}}{1-e^{i(\zeta-\theta)}} \right| \, dk(\theta)\\
&=\frac{1}{2}\int \ln \left(\frac{(1-(1+\epsilon)e^{i(\zeta-\theta)})(1-(1+\epsilon)e^{-i(\zeta-\theta)})}{(1-e^{i(\zeta-\theta)})(1-e^{-i(\zeta-\theta)})}\right)dk(\theta)\\
&= \frac{1}{2} \int \ln \left( 1 + \epsilon + \frac{\epsilon^{2}}{2-2\cos(\zeta-\theta)} \right) \, dk(\theta).
\end{aligned}
$$
Therefore,
$$\begin{aligned}
\tilde{\gamma}((1+\epsilon)e^{i\zeta})&\geq \frac{1}{2}\int_{\zeta-\epsilon}^{\zeta+\epsilon} \ln \left( 1+\epsilon+\frac{\epsilon^{2}}{2-2\cos(\zeta-\theta)} \right) \, dk(\theta)\\
&\geq \frac{1}{2}\ln \left( 1 + \epsilon+\frac{\epsilon^{2}}{2-2\cos\epsilon} \right) k(\zeta-\epsilon,\zeta+\epsilon),
\end{aligned}
$$
take $\epsilon$ sufficiently small, this implies
$$
\tilde{\gamma}((1+\epsilon)e^{i\zeta}) \geq \frac{\ln 2}{2}k(\zeta-\epsilon,\zeta+\epsilon).
$$
Thus, by \eqref{LErel} and \eqref{growth}, we have
$$
k(\zeta-\epsilon,\zeta+\epsilon)\leq \frac{2}{\ln 2}(\gamma((1+\epsilon)e^{i\zeta})+\frac{1}{2}\log (1+\epsilon))\leq C\epsilon^{\frac{1}{2}}
$$
for any $e^{i\zeta} \in \Sigma$, concluding the proof. %Since $k$ is locally constant in $\partial\mathbb{D}\backslash\Sigma$, this means precisely that $k$ is $\frac{1}{2}-$H\"older continuous.
\end{pf}

We also need the following general property of the \textit{Lyapunov exponent}. One can find this result in \cite[Proposition~8.1.10]{Simon1}.

\begin{Lemma}\label{LBE}
Let $e^{i\zeta}\in\partial \mathbb{D}$, then for any $\delta\geq 0$, we have
$$
\tilde{\gamma}((1+ \delta)e^{i\zeta}) \geq \ln (1+\delta).
$$
\end{Lemma}

With this result, we can prove the following:

\begin{Lemma}\label{LB}
For $e^{i\zeta} \in \Sigma$ and sufficiently small $\epsilon > 0$, we have
$$
k(\zeta -\epsilon,\zeta +\epsilon) \geq c \epsilon^{\frac{3}{2}}.
$$
\end{Lemma}

\begin{pf}
Let $\delta = c\epsilon^{\frac{3}{2}}$. For $e^{i\zeta} \in \Sigma$, we have $\tilde{\gamma} (e^{i\zeta}) = 0$ by Lemma~\ref{Prop3.2}. Thus, by a calculation from the proof of Corollary~\ref{upper}, we have
\begin{equation}\label{THappl1}
\tilde{\gamma} ((1+\delta) e^{i\zeta}) =\frac{1}{2} \int \ln \left( 1 + \delta + \frac{\delta^{2}}{2-2\cos(\theta-\zeta)} \right) \, dk(\theta).
\end{equation}
Partition the integration region into $\{ \vert \theta-\zeta\vert \ge \frac{\pi}{10000} \}$, $\{ \epsilon \le \vert \theta-\zeta\vert < \frac{\pi}{10000} \}$, $\{ \epsilon^{4} \le \vert\theta-\zeta\vert < \epsilon \}$, $\{ \vert\theta-\zeta\vert < \epsilon^{4} \}$, so that with the corresponding partition of the RHS of \eqref{THappl1} into the four terms $I_1, I_2, I_3, I_4$, we have $\tilde{\gamma}((1+\delta)e^{i\zeta})=I_{1}+I_{2}+I_{3}+I_{4}$.

As
$$
I_{1} \leq \frac{1}{2} \int \ln \left( 1 + \delta + \frac{\delta^{2}}{2-2\cos \frac{\pi}{10000}} \right) \, dk(\theta),
$$
for $\delta$ sufficiently small, we have $I_{1}\leq \frac{1}{2}\delta$.

Furthermore, we have
$$
\begin{aligned} I_{4} & = \sum_{k\geq 4} \int_{\epsilon^{k} > \vert \theta-\zeta\vert \ge \epsilon^{k+1}} \ln \left( 1 + \delta + \frac{\delta^{2}}{2-2\cos(\theta-\zeta)} \right) \, dk(\theta) \\
& \leq \sum_{k \geq 4} \epsilon^{\frac{k}{2}} \ln \left( 1 + \frac{2\delta^{2}}{\epsilon^{2k+2}} \right) \\ %\, dk(\theta') \\
& \leq \sum_{k \geq 4} \epsilon^{\frac{k}{2}} \ln (1 + 2c^{2} \epsilon^{-2k+1}) \\
& \leq \epsilon^{\frac{7}{4}}.\end{aligned}
$$

With $m=[\ln \frac{\pi\epsilon^{-1}}{10000}]$, we have
$$
I_{2} \leq \sum_{k=0}^{m} \int_{\frac{\pi}{10000} e^{-k-1} \le |\theta-\zeta| < \frac{\pi}{10000} e^{-k}} \ln \left( 1 + \delta + \frac{\delta^{2}}{2-2\cos(\theta-\zeta)} \right) \, dk(\theta').
$$
Suppose $k_{0}>0$ is such that
$$
\delta > \frac{\delta^{2}}{2-2\cos \frac{\pi e^{-k-1}}{10000}} \quad \text{ for } k\leq k_{0},
$$
and
$$
\delta \leq \frac{\delta^{2}}{2-2\cos \frac{\pi e^{-k-1}}{10000}} \quad \text{ for } k_{0} < k \leq m.
$$
Then,
$$
\begin{aligned}I_{2}&\leq 2\delta\sum_{k=0}^{k_{0}}\frac{\pi^{\frac{1}{2}}}{100}e^{-\frac{k}{2}}+2\sum_{k=k_{0}+1}^{m}\frac{\pi^{\frac{1}{2}}}{100}e^{-\frac{k}{2}}\frac{\delta^{2}}{2-2\cos\frac{\pi}{10000}e^{-k-1}}\\
& \leq 2\delta \frac{\pi^{\frac{1}{2}}}{100} \frac{1}{1-e^{-\frac{1}{2}}} + 2\delta^{2} \sum_{k=1}^{m} \frac{\pi^{\frac{1}{2}}}{100} \frac{10000^{2} e^{\frac{3k}{2}+2}}{\pi^{2}} \\
& \leq \frac{\delta}{10} + 2 \times 10^{6} \delta^{2} e^{2} \pi^{-\frac{3}{2}} \frac{1-e^{\frac{3}{2}m}}{1-e^{\frac{3}{2}}} \\
& < \frac{\delta}{10} + C \delta^{2} e^{\frac{3m}{2}} \\
& \leq \frac{\delta}{10} + C c \delta.
\end{aligned}
$$
It follows that
$$
I_{3} \geq \tilde{\gamma}((1+\delta)e^{i\zeta}) - \frac{6\delta}{10} - Cc\delta.
$$
By Lemma~\ref{LBE} we have
$$
\tilde{\gamma}((1+\delta)e^{i\zeta}) \geq \ln (1 + \delta).
$$
Since the constant $c$ above is consistent with our choice of $\delta$, we can therefore shrink it such that $I_{3} \geq \frac{1}{10}\delta$. Since $I_{3} \leq C k(\zeta -\epsilon,\zeta +\epsilon)\ln\epsilon^{-1}$, the result follows.
\end{pf}

\section{Proof of the Main Theorem}\label{sec.5}

Our aim is to prove Theorem~\ref{THM4.3} in the present section. As this theorem makes statements both about the discrete spectrum and the essential spectrum, we divide the discussion into two parts and begin with the easier issue of understanding the discrete spectrum. In fact, everything we need will follow quickly from the general theory. The work from the previous sections then comes into play when we turn our attention to the essential spectrum.

\subsection{Discrete Eigenvalues}

In this subsection we begin with a discussion of the discrete spectrum. As it is well known that for every $x \in \mathbb{T}^d$, $\sigma(\mathcal{E}_x) = \sigma_\mathrm{ess}(\mathcal{E}_x) = \Sigma$, and hence the discrete spectrum of $\mathcal{E}_x$ is empty, we can limit our attention to the discrete spectrum of $\mathcal{C}_x$.

\begin{Theorem}
For every $x \in \mathbb{T}^d$, $\sigma_\mathrm{ess}(\mathcal{C}_x) = \Sigma$, and there is at most one eigenvalue of $\mathcal{C}_{x}$ in each connected component of $\partial\mathbb{D}\backslash\Sigma$.
\end{Theorem}

\begin{pf}
By \cite[Theorem 4.5.2]{Simon1}, $\mathcal{E}_{x} - Q=\mathcal{C}_{x}\bigoplus\mathcal{K}_{x}$, where $Q$ is trace class and $\mathcal{C}_{x},\mathcal{K}_{x}$ are half-line CMV matrices, acting on $\ell^{2}(\mathbb{Z}_{+})$ and $\ell^{2}(\{\cdots,-2,-1\})$, respectively. Thus, $\sigma_\mathrm{ess}(\mathcal{E}_x) = \sigma_{ess}(\mathcal{C}_{x})\cup \sigma_{ess}(\mathcal{K}_{x})$ (cf.~\cite[Theorem 3]{CP}).
Moreover, $\sigma_{ess}(\mathcal{C}_{x}) = \sigma_{ess}(\mathcal{K}_{x})$ as the translation $Tx=x+\omega$ is minimal whenever $\omega$ is rationally independent. Combining this with $\sigma_\mathrm{ess}(\mathcal{E}_x) = \Sigma$, we obtain $\sigma_{ess}(\mathcal{C}_{x}) = \Sigma$.

According  to \cite[Theorem 10.16.3]{Simon2}, in each connected component of $\partial \mathbb{D} \backslash \Sigma$, $\mathcal{C}_{x}$ has at most one eigenvalue. Note that \cite[Theorem 10.16.3]{Simon2} only gives this result for almost every $x\in\mathbb{T}^{d}$, as in the general dynamically defined setting, the condition $\supp(dk) = \sigma(\mathcal{E}_{x})$ only holds for almost every $x\in\mathbb{T}^{d}$. But in our situation we have $\sigma_{ess}(\mathcal{C}_{x}) = \Sigma$ for each $x \in \mathbb{T}^{d}$, and hence via the same line of reasoning, the statement can be derived for each $x\in\mathbb{T}^{d}$.
\end{pf}

This completes the proof of the statements in Theorem~\ref{THM4.3} concerning spectral parameters in $\partial\mathbb{D}\backslash\Sigma$. It remains to prove the pure absolute continuity of all spectral measures, of both $\mathcal{C}_{x}$ and $\mathcal{E}_{x}$ for arbitrary $x \in \mathbb{T}^d$, inside the common essential spectrum $\Sigma$.

\subsection{Bounded Eigenfunctions and Absolutely Continuous Spectrum}

In order to prepare for our proof of the desired absolute continuity statement, which will be given in the next subsection, we discuss in this subsection how to connect the problem to the cocycle norm estimates we have obtained so far.

Given $\phi\in\partial\mathbb{D}$, let $\mu_{\phi}$ be a nontrivial probability measure on $\partial\mathbb{D}$ such that $\alpha_{n}(d\mu_{\phi})=\phi\alpha_{n}(d\mu)$. Let
$$
F(z) = \int \frac{e^{i\theta}+z}{e^{i\theta}-z} \, d\mu
$$
and
$$
F^{\phi}(z) = \int \frac{e^{i\theta}+z}{e^{i\theta}-z} \, d\mu_{\phi}
$$
be the Carath\'eodory functions of $\mu,\mu_{\phi}$, respectively. Then \cite[Theorem~3.2.14]{Simon1} implies
\begin{equation}\label{alek}
F^{\phi} = \frac{(1-\phi)+(1+\phi)F}{(1+\phi)+(1-\phi)F}.
\end{equation}
Let $\varphi^{\phi}_{n}(z),\psi^{\phi}_{n}(z)$ be the orthogonal polynomials of $d\mu_{\phi}$. Then according to \cite[Proposition~3.2.1]{Simon1} we have
$$
\begin{pmatrix} \varphi^{\phi}_{n}(z) \\ \bar{\phi}(\varphi^{\phi}_{n})^{*}(z) \end{pmatrix} = \frac{1}{\sqrt{1-|\alpha_{n}|^{2}}} \begin{pmatrix} z & -\bar{\alpha_{n}} \\ -\alpha_{n} z & 1 \end{pmatrix} \begin{pmatrix} \varphi^{\phi}_{n-1}(z) \\ \bar{\phi}(\varphi^{\phi}_{n-1})^{*}(z) \end{pmatrix}
$$
with the initial condition $\begin{pmatrix}1 \\ \bar{\phi}\end{pmatrix}$. Similarly, $\psi^{\phi}_{n}(z)$ obeys
$$
\begin{pmatrix} \psi^{\phi}_{n}(z) \\ -\bar{\phi}(\psi^{\phi}_{n})^{*}(z) \end{pmatrix} = \frac{1}{\sqrt{1-|\alpha_{n}|^{2}}} \begin{pmatrix} z & -\bar{\alpha_{n}} \\ -\alpha_{n} z & 1 \end{pmatrix} \begin{pmatrix} \psi^{\phi}_{n-1}(z) \\ -\bar{\phi}(\psi^{\phi}_{n-1})^{*}(z) \end{pmatrix}
$$
with initial condition $\begin{pmatrix}1 \\ -\bar{\phi}\end{pmatrix}$.

The following result is a consequence of the Jitomirskaya-Last inequality, which was first proven in the Schr\"odinger case \cite{JL} and then worked out in the OPUC setting \cite{Simon2}. Its purpose is to connect the growth of the cocycle with the Carath\'eodory function of the corresponding spectral measure.

\begin{Proposition}\label{JLeq}
For any $e^{i\zeta} \in \partial \D$ and $0 < \epsilon < 1$, we have the following:
\begin{equation}\label{cocycleBounds}
\sup_{\phi\in\partial\mathbb{D}} \vert F^{\phi}((1-\epsilon)e^{i\zeta})\vert\leq C \sup_{0\leq s\leq c\epsilon^{-1}} \left\| A_{e^{i\zeta}}^{s}\right\|^{2}_{0}.
\end{equation}
\end{Proposition}

Since we did not find this statement in the exact form we need, we include an explanation of how to derive it from the OPUC version of the Jitomirskaya-Last inequality for the convenience of the reader.

For $l \in(0,\infty)$, define
$$
\left\| a\right\|_{l}^{2}=\sum_{j=0}^{[l]}\vert a_{j}\vert^{2}+(l-[l])\vert a_{[l]+1}\vert^{2}.
$$
Let $\varphi^{\phi}(z)=\{\varphi_{0}^{\phi}(z),\varphi_{1}^{\phi}(z),\cdots \},\psi^{\phi}(z)=\{\psi_{0}^{\phi}(z),\psi_{1}^{\phi}(z),\cdots \}$. Then it is known that $\left\| \varphi^{\phi}\right\|_{l}\left\| \psi^{\phi}\right\|_{l}$ runs from $1$ to $\infty$ in a strictly monotone way as $l$ runs from $0$ to $\infty$. The CMV version of the Jitomirskaya-Last inequality now reads as follows (cf.~\cite[Theorem~10.8.2]{Simon2}):

\begin{Lemma}\label{JL1}
For $e^{i\zeta} \in \partial \D$ and $r\in[0,1)$, define $l(r)$ to be the unique solution of
\begin{equation}\label{2.109}
(1-r)\left\| \varphi^{\phi}(e^{i\zeta})\right\|_{l(r)}\left\| \psi^{\phi}(e^{i\zeta})\right\|_{l(r)}=\sqrt{2}.
\end{equation}
Then
\begin{equation}\label{2.110}
A^{-1}\left [\frac{\left\| \psi^{\phi}(e^{i\zeta})\right\|_{l(r)}}{\left\| \varphi^{\phi}(e^{i\zeta})\right\|_{l(r)}}\right]\leq\vert F^{\phi}(re^{i\zeta})\vert\leq A\left[\frac{\left\| \psi^{\phi}(e^{i\zeta})\right\|_{l(r)}}{\left\| \varphi^{\phi}(e^{i\zeta})\right\|_{l(r)}}\right],
\end{equation}
where A is a universal constant.
\end{Lemma}

With the help of this lemma we can now prove Proposition~\ref{JLeq}.

\begin{pf}
Applying Lemma~\ref{JL1} with $r=1-\epsilon$, \eqref{2.109} turns into
$$
\left\| \varphi^{\phi}(e^{i\zeta}) \right\|_{l(\epsilon)}\left\| \psi^{\phi}(e^{i\zeta}) \right\|_{l(\epsilon)} = \frac{\sqrt{2}}{\epsilon},
$$
the pair of inequalities in \eqref{2.110} becomes
\begin{equation}\label{2.111}
(\sqrt{2}A)^{-1}\epsilon\left\| \psi^{\phi}(e^{i\zeta})\right\|_{l(\epsilon)}^{2}\leq\vert F^{\phi}((1-\epsilon)e^{i\zeta})\vert\leq \frac{A}{\sqrt{2}}\epsilon\left\|\psi^{\phi}(e^{i\zeta})\right\|_{l(\epsilon)}^{2}.
\end{equation}
Since
$$
\begin{pmatrix}\varphi^{\phi}_{n}(e^{i\zeta})\\ \bar{\phi}(\varphi^{\phi}_{n}(e^{i\zeta}))^{*}\end{pmatrix}=e^{i\frac{n\zeta}{2}}A^{n}_{e^{i\zeta}}\begin{pmatrix}1 \\ \bar{\phi}\end{pmatrix},\begin{pmatrix}\psi^{\phi}_{n}(e^{i\zeta})\\ -\bar{\phi}(\psi^{\phi}_{n}(e^{i\zeta}))^{*}\end{pmatrix}=e^{i\frac{n\zeta}{2}}A^{n}_{e^{i\zeta}}\begin{pmatrix}1\\ -\bar{\phi}\end{pmatrix},
$$
a direct calculation shows that
\begin{equation}\label{2.112}
\left\| \psi^{\phi}(e^{i\zeta})\right\|_{l(\epsilon)}^{2}\leq C\sum_{j=0}^{j=[l(\epsilon)]+1}\left\| A_{e^{i\zeta}}^{j}\right\|^{2}_{0}.
\end{equation}

We need an explicit upper bound on $l(\epsilon)$. Since
$$
\det{\begin{pmatrix} \varphi^{\phi}_{n}(e^{i\zeta}) &\psi^{\phi}_{n}(e^{i\zeta}) \\ \bar{\phi}(\varphi^{\phi}_{n}(e^{i\zeta}))^{*} &-\bar{\phi}(\psi^{\phi}_{n}(e^{i\zeta}))^{*} \end{pmatrix}} = -\bar{\phi}\varphi^{\phi}_{n}(e^{i\zeta})(\psi^{\phi}_{n}(e^{i\zeta}))^{*}-\bar{\phi}\psi^{\phi}_{n}(e^{i\zeta})(\varphi^{\phi}_{n}(e^{i\zeta}))^{*}
$$
and
$$
\det{\begin{pmatrix} \varphi^{\phi}_{n}(e^{i\zeta}) &\psi^{\phi}_{n}(e^{i\zeta}) \\ \bar{\phi}(\varphi^{\phi}_{n}(e^{i\zeta}))^{*} &-\bar{\phi}(\psi^{\phi}_{n}(e^{i\zeta}))^{*}
\end{pmatrix}} = \det(e^{i\frac{n\zeta}{2}}A_{e^{i\zeta}}^{n} \begin{pmatrix} 1 & 1 \\ \bar{\phi} & -\bar{\phi} \end{pmatrix})=(-2\bar{\phi})e^{in\zeta},
$$
and noting that $(\varphi^{\phi}_{n}(e^{i\zeta}))^{*}=e^{in\zeta}\overline{\varphi^{\phi}_{n}(\frac{1}{e^{-i\zeta}})}=e^{in\zeta}\overline{\varphi^{\phi}_{n}(e^{i\zeta})}$ and $(\psi^{\phi}_{n}(e^{i\zeta}))^{*}=e^{in\zeta}\overline{\psi^{\phi}_{n}(e^{i\zeta})}$, we have
\begin{equation}\label{JL}
\varphi^{\phi}_{n}(e^{i\zeta})\overline{\psi^{\phi}_{n}(e^{i\zeta})}+\psi^{\phi}_{n}(e^{i\zeta})\overline{\varphi^{\phi}_{n}(e^{i\zeta})}=2.
\end{equation}
For a positive integer $L$, define
$$
f^{\phi} = (\varphi^{\phi}_{1}, \overline{\varphi^{\phi}}_{1}, \dots, \varphi^{\phi}_{L}, \overline{\varphi^{\phi}}_{L}), \; g^{\phi} = (\overline{\psi^{\phi}}_{1}, \psi^{\phi}_{1}, \dots, \overline{\psi^{\phi}}_{L},\psi^{\phi}_{L}).
$$
By the Cauchy-Schwarz inequality and \eqref{JL},
$$
\left\| f^{\phi} \right\|_{L} \left\| g^{\phi} \right\|_{L} \geq |\langle f^{\phi}, g^{\phi} \rangle| = \left| \sum_{j=1}^{L} \varphi^{\phi}_{j} \overline{\psi^{\phi}}_{j} + \overline{\varphi^{\phi}}_{j} \psi^{\phi}_{j} \right| = 2L.
$$
Therefore,
$$
\left\| \psi^{\phi}(e^{i\zeta})\right\|_{L}\left\| \varphi^{\phi}(e^{i\zeta})\right\|_{L}\geq L,
$$
which implies that $l(\epsilon)\leq \sqrt{2}\epsilon^{-1}+1$. Thus ($\ref{2.112}$) turns into
$$
\left\| \psi^{\phi}(e^{i\zeta})\right\|_{l(\epsilon)}^{2}\leq C\sum_{0}^{c\epsilon^{-1}}\left\| A_{e^{i\zeta}}^{j}\right\|^{2}_{0},
$$
and the second inequality in \eqref{2.111} becomes
$$
\vert F^{\phi}((1-\epsilon)e^{i\zeta})\vert \leq C \sup_{0 \leq s \leq c\epsilon^{-1}} \left\| A_{e^{i\zeta}}^{s}\right\|^{2}_{0},
$$
concluding the proof.
\end{pf}

Recall from Subsection~\ref{MSM} that for $x \in \mathbb{T}^d$, $\mu_x$ and $\Lambda_x$ denote the canonical maximal spectral measures of $\mathcal{C}_{x}, \mathcal{E}_{x}$, respectively. As a consequence of Proposition~\ref{JLeq}, we have the following:

\begin{Proposition}\label{lm1}
For any $x \in \mathbb{T}^d$ and any $e^{i\zeta} \in \partial\mathbb{D}$ and $0 < \epsilon < 1$, we have
\begin{equation}\label{measure1}
\mu_x(\zeta-\epsilon, \zeta+\epsilon) \leq C\epsilon \sup_{0 \leq s \leq c \epsilon^{-1}} \left\| A^{s}_{e^{i\zeta}} \right\|^{2}_{0}
\end{equation}
and
\begin{equation}\label{measure2}
\Lambda_x(\zeta-\epsilon, \zeta+\epsilon) \leq C \epsilon \sup_{s \leq c \epsilon^{-1}} \left\| A^{s}_{e^{i\zeta}} \right\|_{0}^{2},
\end{equation}
where $C>0$ is a universal constant.
\end{Proposition}

\begin{pf}
As $x$ is fixed throughout the proof, let us drop this parameter from the notation.

\textbf{Half-line case:} Let
$$
F(z) = \int \frac{e^{i\theta}+z}{e^{i\theta}-z} \, d\mu (\theta)
$$
be the Carath\'eodory function of the (canonical maximal) spectral measure $\mu$ of $\mathcal{C}$. Then,
$$
\Re F(re^{i\zeta}) = \int \frac{1-r^{2}}{1+r^{2}-2r\cos(\theta-\zeta)} \, d\mu(\theta).
$$
For $r=1-\epsilon$, we have
$$
\Re F((1-\epsilon) e^{i\zeta}) = \int \frac{2\epsilon-\epsilon^{2}}{2-2\epsilon+\epsilon^{2}-2(1-\epsilon) \cos(\theta-\zeta)} \, d\mu(\theta).
$$
Consequently,
$$
\Re F((1-\epsilon) e^{i\zeta}) \geq \int_{\zeta-\epsilon}^{\zeta+\epsilon} \frac{2\epsilon-\epsilon^{2}}{2 - 2\epsilon + \epsilon^{2} - 2(1-\epsilon)(1 - \frac{\epsilon^{2}}{2})} \, d\mu(\theta),
$$
which implies that
$$
\Re F((1-\epsilon)e^{i\zeta}) \geq \frac{1}{\epsilon}(\mu(\zeta-\epsilon,\zeta+\epsilon)).
$$
By Proposition \ref{JLeq} ($\phi=1$), the result then follows.

\textbf{Full-line case:} For the full-line matrix $\mathcal{E}_{x}$, if we modify $\alpha_{-1}$ into some $\tilde{\alpha}_{-1}\in\partial\mathbb{D}$, then $\mathcal{E}_{x}$ decouples into two half-line matrices $\mathcal{C}^{+}_{x}$ and $\mathcal{C}^{-}_{x}$. Let $F_{+}(z), F_{-}(z)$ be the Carath\'eodory functions of their spectral measure $\mu_{+},\mu_{-}$ respectively. The idea is to apply the CMV version of the Damanik-Killip-Lenz maximum modulus principle argument \cite{DKL} of Munger-Ong \cite{MO}, which uses the Alexandrov family of half-line CMV matrices and the maximum modulus principle to control the Carath\'eodory function of $\mathcal{E}_{x}$.

Define the anti-Carath\'eodory function
$$
M_{-}(z) = \frac{\Re (1-\bar{\alpha}_{0}) - i \Im(1 + \bar{\alpha}_{0}) F_{-}(z)}{i\Im(1 - \bar{\alpha}_{0}) - \Re(1 + \bar{\alpha}_{0})F_{-}(z)}
$$
and let $G(z;k,l)=\langle\delta_{k},(\mathcal{E}_{x}-z)^{-1}\delta_{l}\rangle$. Then,
\begin{equation}\label{der1}
\vert G(z;0,0) + G(z;1,1)\vert \leq \left\vert\frac{1-F_{+}(z)M_{-}(z)}{F_{+}(z)-M_{-}(z)} \right\vert.
\end{equation}
Moreover, the Carath\'eodory function of $\mathcal{E}_{x}$ is
$$
\Phi(z) = \int \frac{e^{i\theta}+z}{e^{i\theta}-z} \, d\Lambda(\theta),
$$
where $\Lambda$ is the canonical maximal spectral measure (cf.~Subsection~\ref{MSM}). Then,
\begin{equation}\label{der2}
\Phi(z) = 1 + 2z (G(z;0,0)+G(z;1,1)).
\end{equation}
Write
$$
\frac{1-F_{+}(z)M_{-}(z)}{F_{+}(z)-M_{-}(z)} = \frac{(1-\frac{M_{-}(z)+1}{M_{-}(z)-1}) + (1 + \frac{M_{-}(z) + 1}{M_{-}(z) - 1}) F_{+}(z)}{(1 + \frac{M_{-}(z) + 1}{M_{-}(z) - 1}) + (1 - \frac{M_{-}(z) + 1}{M_{-}(z) - 1}) F_{+}(z)}.
$$
Since $-M_{-}(z)$ is a Carath\'eodory function, we have $\Re M_{-}(z)<0$, which in turn implies that $\frac{M_{-}(z)+1}{M_{-}(z)-1}\in\mathbb{D}$. By the maximum modulus principle,
$$
\left\vert \frac{1 - F_{+}(z)M_{-}(z)}{F_{+}(z) - M_{-}(z)} \right\vert \leq \sup_{\phi \in \partial \mathbb{D}} \left\vert \frac{(1-\phi) + (1+\phi)F_{+}(z)}{(1+\phi) + (1-\phi)F_{+}(z)} \right\vert.
$$
Recalling \eqref{alek}, we thus have
\begin{equation}\label{twoSidedBound1}
\left\vert \frac{1-F_{+}(z)M_{-}(z)}{F_{+}(z)-M_{-}(z)} \right\vert \leq \sup_{\phi\in\partial\mathbb{D}} \left\vert F_{+}^{\phi}(z) \right\vert,
\end{equation}
where $F_{+}^{\phi}$ is the Carath\'eodory function of the Alexandrov measure $d\mu_{+,\phi}$. Therefore, taking $z=(1-\epsilon)e^{i\zeta}$ and using \eqref{der2}, \eqref{der1}, \eqref{twoSidedBound1}, and \eqref{cocycleBounds}, we find that
$$
\begin{aligned}|\Phi((1-\epsilon)e^{i\zeta})| & \leq 1 + 2|G_{00}((1-\epsilon)e^{i\zeta}) + G_{11}((1-\epsilon)e^{i\zeta})| \\
& \leq 1 + 2 \sup_{\phi \in \partial \mathbb{D}}|F_{+}^{\phi}((1-\epsilon)e^{i\zeta})| \\
& \leq C \sup_{0\leq s\leq c\epsilon^{-1}} \|A^{s}_{e^{i\zeta}}\|^{2} .\end{aligned}
$$
A similar argument gives
$$
\Re \Phi((1-\epsilon)e^{i\zeta}) \geq \frac{1}{\epsilon} \Lambda(\zeta-\epsilon, \zeta+\epsilon).
$$
Together they yield
\begin{equation}\label{twoSidedBound2}
\Lambda(\zeta-\epsilon,\zeta+\epsilon)\leq C \epsilon\sup_{0\leq s\leq c\epsilon^{-1}}\|A_{e^{i\zeta}}^{s}\|^{2},
\end{equation}
finishing the proof.
\end{pf}

\subsection{Absolutely Continuous Spectral Measures}

Since the proof of absolute continuity in the half-line case is very similar to (and even easier than) the proof in the full-line case, we only consider the latter scenario and prove the absolute continuity of $\Lambda$ on $\Sigma$.

Let us first recall the following consequence of subordinacy theory, compare \cite[Theorem~10.9.4]{Simon2} and \cite[Corollary~2.4]{GDO}.

\begin{Theorem}
Define
$$
\mathcal{B} = \left\{ e^{i\zeta} \in \partial\mathbb{D} : \limsup \limits_{s\geq 0} \Vert A_{e^{i\zeta}}^{s} \Vert_{0} < \infty \right\},
$$
then both $\mu|_{\mathcal{B}}$ and $\Lambda|_{\mathcal{B}}$ are absolutely continuous.
\end{Theorem}

We mention that subordinacy theory was originally developed by Gilbert-Pearson in the setting of half-line continuum Schr\"odinger operators \cite{GP} and then extended to the whole-line case by Gilbert \cite{G89}. The relevance of bounded solutions in this context was pointed out by Behnke \cite{B91}, Simon \cite{S96}, and Stolz \cite{S92}.

We also need the following discrete analog of a result of Eliasson from \cite{Eli92}:

\begin{Theorem}\label{Eli92}
Let $\delta > 0$, $\omega \in \mathrm{DC}(\kappa,\tau)$, and $A_{0}\in \mathrm{SL}(2,\mathbb{R})$. Then there is a constant $\epsilon = \epsilon(\gamma,\tau,\delta,\| A_{0}\|)$ such that if $A \in C^{\omega}_{\delta}(\mathbb{T}^{d},\mathrm{SL}(2,\mathbb{R}))$ is real analytic with
$$
\| A - A_{0} \|_{\delta} \leq \epsilon
$$
and the rotation number of the cocycle $(\omega, A)$ satisfies
$$
\|2\rho(\omega,A) - \langle n, \omega \rangle\|_{\mathbb{R}/\mathbb{Z}} \geq \frac{\kappa}{|n|^{\tau}} \quad \forall \;  0 \neq n \in \mathbb{Z}^{d}
$$
or $2\rho(\omega,A) - \langle n, \omega \rangle\in\mathbb{Q}$ for some $n \in \mathbb{Z}^{d}$, then $(\omega,A)$ is reducible to constant coefficients of a quasi-periodic (perhaps with frequency $\frac{\omega}{2}$) and analytic transformation.
\end{Theorem}

It is enough to prove that
\begin{equation}\label{e.sigmabgoal}
\Lambda (\Sigma \backslash \mathcal{B}) = 0.
\end{equation}
Let $\mathcal{R}$ be the set of $e^{i\zeta} \in \Sigma$ such that the cocycle is reducible. We know that $\mathcal{R} \backslash \mathcal{B}$ only contains elements $e^{i\zeta}$ for which $(\omega,S(\alpha,e^{i\zeta}))$ is analytically reducible to a constant parabolic cocycle. By Theorem~\ref{Eli92} it follows that $\mathcal{R}\backslash\mathcal{B}$ is countable: indeed for any such $e^{i\zeta}$, the well known gap labeling theorem ensures that there exists a $k\in \Z^{d}$ such that $2\rho(\omega,S(\alpha,e^{i\zeta}))=\langle k,\omega\rangle \mod \mathbb{Z}$. If $e^{i\zeta}\in\mathcal{R}$, by \eqref{equivalentmodel}, any non-zero solution of $\mathcal{E}_{x} u = e^{i\zeta}u$ satisfies $\inf \limits_{n\in\Z} \vert u_{n} \vert^{2} + \vert u_{n+1} \vert^{2} > 0$. In particular there are no eigenvalues in $\mathcal{R}$, thus (using countability of $\mathcal{R}\backslash\mathcal{B}$) $\Lambda(\mathcal{R} \backslash \mathcal{B})=0$. It therefore suffices to show
\begin{equation}\label{e.lambdargoal}
\Lambda (\Sigma \backslash \mathcal{R}) = 0.
\end{equation}

Let $J_{j}(e^{i\zeta})$ be an open $2^{\frac{2}{3}} \epsilon_{j-1}^{\frac{2}{45}}$-neighborhood of $e^{i\zeta} \in K_{j}$ in $\partial \D$ (recall that $K_j$ was defined in \eqref{e.Kjdef}). By Corollary~\ref{growth01} and Proposition~\ref{lm1} we have
\begin{align*}
\Lambda(J_{j}(e^{i\zeta})) & \leq \sup_{0\leq s\leq C\epsilon_{j-1}^{-\frac{2}{45}}}\|A_{e^{i\zeta}}^{s}\|^{2}_{0}\vert J_{j}(e^{i\zeta})\vert \\
& \leq \sup_{0\leq s\leq C\epsilon_{j-1}^{-\frac{1}{16}}}\|A_{e^{i\zeta}}^{s}\|^{2}_{0}\vert J_{j}(e^{i\zeta})\vert \\
& \leq C\epsilon_{j-1}^{\frac{1}{48}},
\end{align*}
where $|\cdot|$ denotes Lebesgue measure. Take a finite subcover such that $\overline{K_{j}}\subset \bigcup_{l=0}^{r}J_{i}(e^{i\zeta_{l}})$. Refining this subcover if necessary, we may assume that every $z \in \partial \D$ is contained in at most 2 different $J_{j}(e^{i\zeta_{l}})$. By Lemma~\ref{LB} and \eqref{rotationDOS},
$$
\vert 2\rho(J_{j}(e^{i\zeta})) \vert = k(\zeta-2^{-\frac{1}{3}}\epsilon_{j-1}^{\frac{2}{45}}, \zeta+2^{-\frac{1}{3}}\epsilon_{j-1}^{\frac{2}{45}}) \geq 2 c \epsilon_{j-1}^{\frac{1}{15}}.
$$
By Lemma~\ref{LM4}, if $e^{i\zeta}\in K_{j}$, we have
$$
\left\| 2\rho(e^{i\zeta})-\langle n_{j},\omega\rangle\right\|_{\R/\Z}\leq 2\epsilon_{j-1}^{\frac{1}{15}}
$$
for some $\vert n_{j}\vert <2N_{j-1}$. This implies that $2\rho(K_{j})$ can be covered by $2N_{j-1}$ open arcs $T_{s}$ of length $2\epsilon_{j-1}^{\frac{1}{15}}$. Since $\vert T_{s}\vert \leq \frac{1}{c}\vert 2\rho(J_{j}(e^{i\zeta}))\vert$ for any $s$, $e^{i\zeta}\in K_{j}$, there are at most $2([\frac{1}{c}]+1)+4$ open arcs $J_{j}(e^{i\zeta_{l}})$ such that $2\rho(J_{j}(e^{i\zeta_{l}}))$ intersects $T_{s}$. We conclude that there are at most $2(2([\frac{1}{c}]+1)+4)N_{j-1}$ open arcs $J_{j}(z_{l})$ to cover $K_{j}$. Then
\begin{equation}\label{boundOfMeas}
\Lambda(K_{j})\leq \sum_{j=0}^{r}\Lambda(J_{j}(e^{i\zeta_{l}})) \leq CN_{j-1}\epsilon_{j-1}^{\frac{1}{48}} \leq C\epsilon_{j-1}^{\frac{7}{384}}.
\end{equation}
Since $\epsilon_{j}=\epsilon_{0}^{2^{j}}$ and $\epsilon_{0}$ is small, \eqref{boundOfMeas} implies
\begin{equation}\label{e.BCprep}
\sum\limits_{j} \Lambda(\overline{K_{j}}) < \infty.
\end{equation}
Since $\Sigma \backslash \mathcal{R} \subseteq \limsup K_{j}$, the Borel-Cantelli Lemma and \eqref{e.BCprep} imply \eqref{e.lambdargoal}. By our earlier discussion, this in turn implies \eqref{e.sigmabgoal}. Therefore, $\Lambda = \Lambda|_{\mathcal{B}}$, which is purely absolutely continuous.

This completes the proof of absolute continuity and, together with our discussion above of discrete eigenvalues in the half-line case, the proof of Theorem~\ref{THM4.3}.

\begin{appendix}

\section{A Quantitative Almost Reducibility Result}\label{sec.app}

The following quantitative almost reducibility result from \cite{CCYZ,LYZZ} is the basis of our proof.

\begin{Proposition}\label{rd}
Let $\omega \in \mathrm{DC}(\kappa,\tau)$, $\kappa, \tau, r > 0$, $\sigma = \frac{1}{15}$, $S_{0} \in \mathrm{SU}(1,1)$, $f_{0} \in C^{\omega}_{r}(\T^{d},\mathrm{su}(1,1))$. Then for any $r'\in(0,r)$, there exist a constant $D_{0} = D_{0}(\kappa,\tau)$ and a numerical constant $C_{0}$ such that if
\begin{equation}\label{initialCondition}
\Vert f_{0} \Vert_{r} \leq \epsilon \leq \frac{D_{0}}{\left\| S_{0} \right\|^{C_{0}}} \left( \min \left\{ 1, \frac{1}{r} \right\} (r-r') \right)^{C_{0}\tau},
\end{equation}
then there exists $B\in C^{\omega}_{r'}(2\T^{d},\mathrm{SU}(1,1)),S_{+}\in \mathrm{SU}(1,1),f_{+}\in C^{\omega}_{r'}(\T^{d},\mathrm{su}(1,1))$ such that
$$
B(x+\omega)(S_{0}e^{f_{0}(x)})B^{-1}(x)=S_{+}e^{f_{+}(x)}.
$$
Define $N=\frac{2\vert \ln \epsilon\vert}{r-r'}$, let $spec(S_{0})=\{e^{ 2\pi i\rho},e^{- 2\pi i\rho}\}$, then we have the following estimates:

$\bullet$ \textbf{(Non-resonant case):} If $\vert 2\rho-\langle n,\omega\rangle\vert\geq\epsilon^{\sigma}$ holds for any $ n\in \Z^{d}$ with $0<\vert n\vert <N$, then
$$
\left\| B\right\|_{r'}\leq 1+\epsilon^{\frac{1}{2}},\left\| B\right\|_{0}\leq 1+\epsilon^{\frac{1}{2}},
$$
$$
\|f_{+}\|_{0},\left\| f_{+}\right\|_{r'}\leq e^{-N(r-r')} \leq\epsilon^{2},
$$
$$
\left\| S_{+}-S_{0}\right\|\leq 2\epsilon
$$

$\bullet$ \textbf{(Resonant case):} If $\vert 2\rho-\langle n_{*},\omega\rangle\vert<\epsilon^{\sigma}$ holds for some $n_{*}\in \Z^{d}$ with $0<\vert n_{*}\vert <N$, then
$$
\left\|B\right\|_{r'}\leq C\vert n_{*}\vert^{\tau}\epsilon^{\frac{r'}{r'-r}}, \left\| B\right\|_{0}\leq C\vert n_{*}\vert^{\tau},
$$
$$
\|f_{+}\|_{0},\left\|f_{+}\right\|_{r'}\leq \epsilon e^{-r'\epsilon^{-\frac{1}{18\tau}}}.
$$
Moreover, $\deg B=n_{*}$ and the constant $S_{+}$ can be written as
$$
S_{+}=\exp\begin{pmatrix} i t_{+} & v_{+} \\ \bar{v}_{+} & -i t_{+} \end{pmatrix},
$$
where $t_{+}\in\mathbb{R}, v_{+}\in\mathbb{C}$ with $\vert t_{+}\vert\leq \epsilon^{\frac{1}{16}}, \vert v_{+}\vert\leq \epsilon^{\frac{15}{16}}e^{-2\pi\vert n_{*}\vert r}$.
\end{Proposition}

\section{An Approach to Establishing a Non-Perturbative Result in the Case $d = 1$}\label{app.NP}

In this appendix, we explain how to establish purely absolutely continuous spectrum for $\mathcal{E}_x$ in the non-perturbative regime when the underlying torus dimension is given by $d = 1$. That is, one wishes to strengthen the statement in Theorem \ref{THM4.3} in such a way that the smallness of $\lambda$ does not depend on the Diophantine constants of the frequency $\omega$. The approach we explain here was first developed in \cite{XYZ}, we just give a sketch here for completeness.

First, we need the following:

\begin{Theorem}\label{np}{\rm \cite{HY,YZ}}
Let $\omega \in \R \backslash \Q$, $S_{0} \in \mathrm{SU}(1,1)$, $f_{0} \in C^{\omega}_{r}(\T,\mathrm{su}(1,1))$. There exists $\tilde{\epsilon} = \tilde{\epsilon}(S_0, r) > 0$, such that if $\Vert f_{0} \Vert_{r} \leq \tilde{\epsilon}$, then $(\omega, S_0 e^{f_0})$ is almost reducible.
\end{Theorem}

Therefore, let us consider a Szeg\H{o} cocycle $(\omega, S_{0}(e^{i\zeta})e^{f_{0}(x,\zeta)}) = (\omega, S(\alpha,e^{i\zeta}))$ that is close to constant. If $\lambda$ is small, then one can apply Theorem~\ref{np}, and there exists $\Phi^k_{e^{ i \zeta}} \in  C^{\omega}_{r_{k}}(2\T, \mathrm{SU}(1,1))$ such that
$$
\Phi^k_{e^{i \zeta}}(x+\omega) S_{0}(e^{i\zeta})e^{f_{0}(x,\zeta)} \Phi^k_{e^{ i \zeta}}(x)^{-1} = S_k(e^{i\zeta}) e^{f_k(x,\zeta)}.
$$
If $k$ is large enough, then one can reduce the initial cocycle $(\omega, S_{0}(e^{i\zeta})e^{f_{0}(x,\zeta)})$ to the perturbative regime defined in Proposition~\ref{rd}. Moreover, since the spectrum $\Sigma$ is compact, one can apply the compactness argument from \cite[Proposition 5.2]{LYZZ} and show that there exists $\Gamma = \Gamma(\alpha, r)$, which is independent of $\zeta$, such that
\begin{align}
\|\Phi^k_{e^{i \zeta}}\|_{r_k} & \leq \Gamma, \label{es1}\\
| \deg \Phi^k_{e^{i \zeta}}| & \leq C |\ln \Gamma|. \label{es2}
\end{align}
In other words, after a finite number (that is uniform with respect to $e^{i \zeta} \in \Sigma$) of conjugation steps,   one can reduce the cocycle to the perturbative regime.

Consequently, one can apply Proposition \ref{rd1}, and control the growth of the cocycles in the resonant set $K_j$.  Corollary \ref{growth01} will follow from \eqref{es1}, while Lemma~\ref{LM4} still follows as a result of \eqref{es2}. The rest of the proof will follow the same way as in the main body of the present paper.

We point out that one cannot instead use Eliasson's perturbative result from \cite{Eli92}, since he uses a parameterized KAM method, and there one cannot suitably control the $\zeta$-dependence of $S_k(e^{i\zeta})$.

\section{Quasi-Periodic CMV Matrices With Singular Continuous Spectrum}\label{app.sc}

Let us give a brief discussion of the phenomenon of singular continuous spectrum in the context of analytic\footnote{In the $C^0$ category it is known that generically, the spectrum of an extended quasi-periodic CMV matrix is purely singular continuous: generic absence of point spectrum was shown in \cite{O13} and generic absence of absolutely continuous spectrum was shown in \cite{FDG}.} quasi-periodic extended\footnote{The corresponding problem for standard analytic quasi-periodic CMV matrices, while interesting, is not well understood. The proofs of the results in the extended case do not carry over to the standard case.} CMV matrices. We rely on known results and hence the main purpose is to make explicit how those known results are relevant to this question.

Given the well-known parallels between the theory of discrete one-dimen\-sional Schr\"odinger operators (and more generally Jacobi matrices) and CMV matrices, we are naturally guided by what is known in the Schr\"odinger case. Given this perspective, one should single out two mechanisms that produce examples with singular continuous spectrum in the Schr\"odinger context:
\begin{itemize}

\item[(i)] the coexistence of positive Lyapunov exponents and Liouville frequencies,

\item[(ii)] the self-duality with respect to Aubry duality.

\end{itemize}

In approach (i) one needs a result that establishes positive Lyapunov exponents on the spectrum (or an energy region) for all minimal translations. This input ensures purely singular spectrum. In addition one needs to chose the translation to be of sufficient Liouville nature so that one can exclude eigenvalues via the Gordon lemma \cite{ayz,G}. Let us for simplicity consider the case $d =1 $. Here, to measure how exponentially Liouvillean $\omega$ is, we consider
$$
\beta(\omega):=\limsup_{n\rightarrow \infty}\frac{\ln q_{n+1}}{q_n},
$$
where $\frac{p_n}{q_n}$ are the continued fraction approximants of $\omega$. This input ensures continuous spectrum. As a result one obtains singular continuous spectrum, and it should be emphasized that due to the presence of positive Lyapunov exponents, the spectral measures are highly singular -- they are zero-dimensional and indeed supported by sets of zero capacity \cite{S07}.

In approach (ii) one uses that duality transforms pure point spectrum into absolutely continuous spectrum, and in some weaker sense vice versa as well. Self-duality therefore implies the absence of both. In addition, duality also transforms positive Lyapunov exponents into zero Lyapunov exponents, and hence the singular continuous spectral measures one obtains in this approach are accompanied by zero Lyapunov exponents due to self-duality \cite{AJM,J}. Indeed, they are (expected to be) more regular, with higher dimensionality and (hence the necessity for) larger supports.

\medskip

Turning our attention to the CMV case now, we note that approach (i) can be implemented as there are CMV versions of the positive Lyapunov exponent result \cite{Zhang1} and the Gordon lemma \cite{O12}. Indeed, let us attempt to push the Gordon lemma aspect to its limit. Based on the sharp Gordon lemma from \cite{ayz}, we have the following:

\begin{Theorem}\label{gordon}
Suppose that $\alpha \in C^1(\T, \R)$, $\omega \in \R \backslash \Q$ with $\beta(\omega) > 0$. Then for every $x \in \mathbb{T}$, $\mathcal{E}_x$ has purely singular continuous spectrum on the set
$$
\mathcal{S} = \{ e^{i\zeta} \in \Sigma : \beta(\omega) > \gamma(e^{i\zeta}) > 0 \}.
$$
\end{Theorem}

\begin{pf}
This is a Gordon-type statement that is essentially contained in the proof of \cite[Theorem 1.1]{ayz}. We just give the short argument here for completeness as we are dealing with the extended CMV case (whereas \cite{ayz} considers the Schr\"odinger case).

Given $e^{i\zeta} \in \mathcal{S}$, we have, by a uniform upper semi-continuity and telescoping argument (see, e.g., \cite[Proposition~3.1]{ayz}), that for any sufficiently small $\epsilon > 0$, there exists $K = K(\zeta, \omega, \alpha, \epsilon)$, which is in particular independent of $x$, such that for $n \geq K$, we have
\begin{eqnarray}
\label{appro-3}\sup_{x\in\T}\| A^{q_n}_{e^{i \zeta}}(x+q_n\omega)-
A^{q_n}_{e^{i \zeta}} (x)\|& \leq&
e^{-(\beta-\gamma- \epsilon)q_n}.\\
\label{appro-4}\sup_{x\in\T}\| A^{-q_n}_{e^{i \zeta}}(x+q_n\omega)-
A^{-q_n}_{e^{i \zeta}}(x)\|& \leq& e^{-(\beta-\gamma- \epsilon)q_n}.
\end{eqnarray}
As a consequence, we can apply the following lemma:

\begin{Lemma}\label{l.Gordonlemma}{\rm \cite{AJZ16,ayz}}
Suppose that \eqref{appro-3} and \eqref{appro-4} hold. Then we have
$$%\begin{equation}\label{refinegordon}
\max \left\{ \Big\| A^{q_{n}}_{e^{i \zeta}}(x) \begin{pmatrix} \varphi_{0} \\ \varphi^{*}_{0}\end{pmatrix} \Big\|, \Big\| A^{-q_{n}}_{e^{i \zeta}}(x) \begin{pmatrix} \varphi_{0} \\ \varphi^{*}_{0} \end{pmatrix} \Big\|, \Big\| A^{2q_{n}}_{e^{i \zeta}}(\theta) \begin{pmatrix} \varphi_{0} \\ \varphi^{*}_{0} \end{pmatrix} \Big\| \right\} \geq \frac{1}{2\sqrt{2}},
$$%\end{equation}
where $\begin{pmatrix} \varphi_{0} \\ \varphi^{*}_{0} \end{pmatrix} = \begin{pmatrix} \varphi_{0}(e^{i\zeta}) \\ \varphi^{*}_{0}(e^{i\zeta}) \end{pmatrix}$ is the initial value at $e^{i\zeta}$.
\end{Lemma}

According to \eqref{e.equinorm}, \eqref{e.szegorecnorm}, \eqref{e.Sz}, and noting that
$$
A_{e^{i\zeta}}^{n} = e^{-\frac{i n\zeta}{2}} \tilde{S}(\alpha_{n},e^{i\zeta}) \cdots \tilde{S}(\alpha_{1},e^{i\zeta}),
$$
we are able to utilize the transfer matrix $A_{e^{i\zeta}}^{n}$ to evaluate the norm of  $\begin{pmatrix}s_{n}\\ t_{n}\end{pmatrix}$, the solutions of generalized eigenvalue equations $\mathcal{E}_{x} s = e^{i\zeta} s$. Thus Lemma~\ref{l.Gordonlemma} ensures that for any $x \in \mathbb{T}$, $\mathcal{S}$ does not contain any eigenvalues of $\mathcal{E}_{x}$. The result now follows from the Ishii-Pastur Theorem (cf.~\cite[Theorem 10.5.7]{Simon2}).
\end{pf}

\begin{Remark}
(a) Theorem~\ref{gordon} is the CMV analog of a Schr\"odinger result that is known to be optimal: for the almost Mathieu operator, if $\gamma > \beta(\omega)$, then the operator exhibits Anderson localization \cite{ayz}. One may expect that Theorem~\ref{gordon} is optimal as well.

(b) We emphasize that the result holds for every phase $x \in \mathbb{T}$ and that the proof exploits the two-sided nature of the problem. In view of the latter aspect, we regard it as an interesting problem to prove a singular continuity result for standard (one-sided) CMV matrices with analytic quasi-periodic Verblunsky coefficients.\footnote{Of course we mean that the desired singular continuity of spectral measures should be established on the \textit{essential} spectrum of $\mathcal{C}_x$.}

(c) To make Theorem~\ref{gordon} meaningful, one can apply it, for example, to the specific model studied by Zhang in \cite{Zhang1}, as he established the positivity of the Lyapunov exponent there and $\beta(\omega)$ can be made as large as needed (even infinite) for suitable choices of $\omega$ (that will always form a dense $G_\delta$ set).
\end{Remark}

\medskip

Approach (ii), establishing purely singular continuous spectrum for self-dual models, has been worked out in the unitary case by Fillman-Ong-Zhang \cite{FOZ}, albeit only for almost every phase. It is curious that while the general duality correspondence between the pp and ac regimes is still absent in the extended CMV setting (as pointed out repeatedly in this paper), the somewhat degenerate situation of a self-dual model has been identified in the unitary setting (for quantum walks, which are closely connected to extended CMV matrices via the CGMV connection \cite{CGMV}).

\end{appendix}

\end{document}